\documentclass[12pt,reqno]{amsart}

\usepackage[all]{xy} 

\usepackage{amsmath}
\usepackage{latexsym,amssymb,eucal}
\usepackage{mathrsfs}
\usepackage{amscd}
\usepackage{pb-diagram}
\usepackage{enumerate}
\usepackage{bm}
\allowdisplaybreaks[1]

\usepackage{mathtools}

\usepackage[dvips]{graphics}

\usepackage{multirow,bigdelim}

\newtheorem{lemma}{Lemma}[section]
\newtheorem{theorem}[lemma]{Theorem}
\newtheorem{proposition}[lemma]{Proposition}
\newtheorem{corollary}[lemma]{Corollary}
\newtheorem{question}[lemma]{Question}
\theoremstyle{definition}
\newtheorem{definition}[lemma]{Definition}
\newtheorem{example}[lemma]{Example}
\newtheorem{remark}[lemma]{Remark}

\numberwithin{equation}{section}

\newcommand{\bdf}{\begin{definition}}
\newcommand{\edf}{\end{definition}}
\newcommand{\blem}{\begin{lemma}}
\newcommand{\elem}{\end{lemma}}
\newcommand{\bthm}{\begin{theorem}}
\newcommand{\ethm}{\end{theorem}}
\newcommand{\bpf}{\begin{proof}}
\newcommand{\epf}{\end{proof}}
\newcommand{\bprop}{\begin{proposition}}
\newcommand{\eprop}{\end{proposition}}
\newcommand{\bcor}{\begin{corollary}}
\newcommand{\ecor}{\end{corollary}}
\newcommand{\brem}{\begin{remark}}
\newcommand{\erem}{\end{remark}}
\newcommand{\bquest}{\begin{question}}
\newcommand{\equest}{\end{question}}
\newcommand{\bex}{\begin{example}}
\newcommand{\eex}{\end{example}}

\newcommand{\benu}{\begin{enumerate}\renewcommand{\labelenumi}{{\rm (\arabic{enumi})}}\renewcommand{\itemsep}{0pt}}
\newcommand{\eenu}{\end{enumerate}}

\newcommand{\N}{\mathbb{N}}
\newcommand{\Z}{\mathbb{Z}}
\newcommand{\Q}{\mathbb{Q}}
\newcommand{\R}{\mathbb{R}}
\newcommand{\C}{\mathbb{C}}

\newcommand{\bM}{\mathbb{M}}

\newcommand{\cU}{\mathcal{U}}

\newcommand{\cW}{\mathcal{W}}
\newcommand{\cZ}{\mathcal{Z}}

\newcommand{\e}{\varepsilon}
\newcommand{\dis}{\displaystyle}

\DeclareMathOperator{\id}{id}
\DeclareMathOperator{\tr}{tr}
\DeclareMathOperator{\Tr}{Tr}
\DeclareMathOperator{\real}{Re}
\DeclareMathOperator{\image}{Im}
\DeclareMathOperator{\Ad}{Ad}

\DeclareMathOperator{\Conv}{conv}
\DeclareMathOperator{\Sp}{sp}

\newcommand{\ip}[1]{\mathopen{\langle}#1\mathclose{\rangle}}

\begin{document}
\title{A note on injective factors with trivial bicentralizer}
\author{Rui OKAYASU}
\address{Department of Mathematics Education, Osaka Kyoiku University, Kashiwara, Osaka 582-8582, JAPAN}
\email{rui@cc.osaka-kyoiku.ac.jp}
\date{\today}
\thanks{The author is partially supported by JSPS KAKENHI Grant Number 17K05278.}

\maketitle

\begin{abstract}
We give an alternative proof that an injective factor on a Hilbert space
with trivial bicentralizer is ITPFI.
Our proof is given in parallel with each type of factors
and it is based on the strategy of Haagerup.
As a consequence,
the uniqueness theorem of injective factors
except type III$_0$ follows from Araki-Woods' result.
\end{abstract}


\section{Introduction}

In Connes' fundamental works in operator algebras, 
one of them is proved that injective factors on a separable Hilbert space are hyperfinite.
Another result is that injective factors of type II$_1$, II$_\infty$ and III$_\lambda$ ($\lambda\ne 1$) are classified in \cite{c3}.
The remaining problem of the uniqueness of injective type III$_1$ factor
is solved by Haagerup in \cite{ha3}
by proving the so-called bicentralizer problem in \cite{c4}.

In \cite{ha1}, Haagerup also give another proof of the first result mentioned the above without the automorphism group machinery of Connes.
In his proof, semidiscreteness rather than injectivity is applied.
In \cite{po}, Popa gives the third approach of this result in the case of type II$_1$.

Alternative proofs of the uniqueness of injective factor of type II and type III are also given by Haagerup in \cite{ha1, ha2, ha4}.
For each case, similar techniques are applied. 
The one of important notions is trivial bicentralizer.
It is essential in the case of type III$_1$,
but the equivalent condition of trivial bicentralizer is more important than the original definition in \cite{ha4}.
This condition is similar to Dixmier property.
This property is applied in the case of type II$_1$ \cite{ha1},
and its relative version is applied in the case of type III$_\lambda$ \cite{ha2}.  
The other of important notions is almost unitary equivalence in a Hilbert bimodule.

Moreover, in comparison with his papers \cite{ha1, ha2, ha4},
the uniqueness of injective type II$_1$ factor follows
from Marray-von Neumann's fundamental result of the uniqueness of hyperfinite type II$_1$ factor.
In the case of type III$_\lambda$, it is directly proved that 
an injective type III$_\lambda$ factor is isomorphic to the Powers factor $R_\lambda$.
In the case of type III$_1$, by using Connes-Woods' characterization of ITPFI factors in \cite{cw}, it is proved that an injective type III$_1$ factor is ITPFI. 
Therefore the uniqueness of injective type III$_1$ factor follows from Araki-Woods' result in \cite{aw}.
These are proved by similar arguments, but they are dependent
on the choice of type of a given injective factor.

In this note, 
we give an alternative proof that an injective factor on a Hilbert space
with trivial bicentralizer is ITPFI.
Our proof is given in parallel with each type of factors
and it is based on the strategy of Haagerup.
One of our purposes is to unify his proof.
Here, we remark that there exists an injective type III$_0$ factor,
which is not ITPFI. However the assumption of trivial bicentralizer in the above assertion excludes the case of type III$_0$ from consideration. Namely, every type III$_0$ factor
has non-trivial bicentralizer. This fact may be a folklore among specialists, but we do not find it in literature.
Hence we also give its proof in this note.
To do so, we define the bicentralizer for a general weight
by using the free ultrafilter.
This is inspired by Houdayer-Isono's paper \cite{hi}.
The starting point is 
the semidiscreteness, which is equivalent to 
the injectivity by the work of Wassermann in \cite{was}.
To achieve the above assert,
we need to generalize Haagerup's works.
The one is to obtain an approximate factorization related to the modular automorphism from the semidiscreteness. 
In the case of type III, it relies on the uniqueness of injective type II$_1$ factor in \cite{ha2, ha4}. However, we independently give such approximate factorization for arbitrary injective von Neumann algebra
by combining a number of technics in \cite{ha1}.
The other is almost unitary equivalence in Hilbert bimodules
established in \cite{ha2},
which is a generalization of \cite{ha1}.
However, as in the case of \cite{ha4},
we require a further generalization.
Finally, based on Haagerup's approach,
we give a proof of the above assertion
by using Connes-Woods' characterization of ITPFI factors.


\section{Semidiscreteness with the modular automorphisms}


Let $M$ be a von Neumann algebra.
We denote by $\cU(M)$ the unitary group of $M$.
For a fn (faithful normal) state $\varphi$,
we denote by $\Delta_\varphi$ (resp. $J_\varphi$)
the modular operator (resp. the modular conjugation operator)
associated with $\varphi$, respectively.
We put
\[
\|x\|_\varphi\coloneqq \varphi(x^*x)^{1/2}
\quad
\text{and}
\quad
\|x\|_\varphi^\sharp\coloneqq \varphi\left(\frac{x^*x+xx^*}{2}\right)^{1/2}
\quad
\text{for}\ x\in M.
\]
We denote by $L^2(M, \varphi)$ the standard form for $M$
with the cyclic unit vector $\xi_\varphi$, 
which becomes a normal $M$-$M$ bimodule, 
where the left and right actions are given by
\[
a\xi x\coloneqq aJ_\varphi x^* J_\varphi\xi
\quad\text{for}\ a, x\in M\ \text{and}\ \xi\in L^2(M, \varphi).
\]
The centralizer of $\varphi$ is denoted by $M_\varphi$.
For $m\in\N$, we denote by
$\tr_m$ the normalized tracial state 
and by $\Tr_m$ the canonical trace
with $\Tr_m(1)=m$
on the $m\times m$ matrix algebra $\bM_m$.

If $M$ is semidiscrete,
then the identity map $\id_M$ on $M$ has 
an approximate factorization through 
matrix algebras $\bM_{m(\lambda)}$,
\[
\xymatrix{
 M \ar[rr]^{\id_M} \ar@{-->}[dr]_{S_{\lambda}} &\ar@{}[d]|{\circlearrowright}&M \\
&\bM_{m(\lambda)} \ar@{-->}[ur]_{T_{\lambda}} &
}
\]
where $(S_\lambda)$ and
$(T_\lambda)$ are nets
of ucp (unital completely positive) maps.
The purpose of this section is to show that 
for a given fn state $\varphi$ on $M$
and a positive number $\delta>0$,
one can choose an approximation factorization
such that moreover
\[
\varphi\circ T_\lambda=\psi_\lambda,
\quad
\psi_\lambda\circ S_\lambda=\varphi,
\]
and
\[
\|\sigma_t^\varphi\circ T_\lambda-T_\lambda\circ\sigma_t^{\psi_\lambda}\|
\leq\delta|t|
\quad\text{for}\
t\in\R,
\]
where $(\psi_\lambda)$ is a net of fn states on $\bM_{m(\lambda)}$.

To do so, we will prepare some lemmas, 
which are essentially proved in \cite{ha1}, \cite{ha2} and \cite{ha4}.
The first lemma is given in \cite{ha1}.

\blem[{\cite[Lemma 3.1]{ha1}}]\label{lem:3.1}
Let $\tau$ be a tracial state on $M$,
and $T\colon\bM_m\to M$ be a faithful ucp map.
Put   
\[
\psi(x)\coloneqq \tau\circ T(x)
\quad\text{for}\
x\in\bM_m,
\]
and let $h\in\bM_m$ be the positive element for which
\[
\psi(x) 
=\tr_m(xh)
\quad\text{for}\
x\in\bM_m.
\]
Then
there exists a unique ucp map $S\colon M\to\bM_m$
such that
$\psi\circ S=\tau$ and
\[
\tr_m(x^*h^{1/2}S(y)h^{1/2})
=
\tau(T(x)^*y)
\quad\text{for}\
x\in\bM_m, y\in M.
\]
Moreover 
\[
\|T(x)\|_\tau^2
\leq
\tr_m(x^*h^{1/2}xh^{1/2})
\quad\text{for}\
x\in\bM_m.
\]
\elem

The second lemma is nothing but \cite[Lemma 3.2]{ha1}.
However, we necessarily sketch a proof
to use it in Remark \ref{rem:commute}.

\blem[{\cite[Lemma 3.2]{ha1}}]\label{lem:3.2}
Let $\psi$ be a state on $\bM_m$
of the form
\[
\psi(x)\coloneqq \tr_m(xh)
\quad\text{for}\quad
x\in\bM_m,
\]
where $h$ has strictly positive rational eigenvalues.
Then there exist ucp maps
$S\colon\bM_m\to\bM_p $ and
$T\colon\bM_p\to\bM_m$
such that
$\tr_p\circ S=\psi$,
$\psi\circ T=\tr_p$
and
\[
\|x-T\circ S(x)\|_\psi^\sharp\leq
\|h^{1/2}x-xh^{1/2}\|_2
\quad\text{for}\
x\in\bM_m.
\]
\elem

\bpf
We may assume that $h$ is an 
$m\times m$ diagonal matrix
with strictly positive rational eigenvalues
$\lambda_1,\dots,\lambda_m$.
Choose positive integers $p_1,\dots,p_m$, $p$
such that
\[
\frac{\lambda_i}{m}=\frac{p_i}{p}
\quad\text{for}\
1\leq i\leq m.
\]
Note that
\[
\sum_{i=1}^mp_i=p.
\]
For $1\leq i,j\leq m$,
we define the $p_i\times p_j$ matrix $F_{ij}$ by
\[
[F_{ij}]_{kl}\coloneqq \delta_{kl}
\quad\text{for}\
1\leq k\leq p_i,
1\leq l\leq p_j,
\]
and the $p\times p$ matrix $f_{ij}$
with block matrix by 
\[
[f_{ij}]_{kl}=\delta_{ik}\delta_{jl}F_{ij}
\quad\text{for}\
1\leq k,l\leq m.
\]
Let $(e_{ij})$ be the matrix units for $\bM_m$.
We define the ucp map 
$S\colon\bM_m\to\bM_p$
by
\[
S\left(
\sum_{i,j=1}^mx_{ij}e_{ij}
\right)
\coloneqq 
\sum_{i,j=1}^mx_{ij}f_{ij}.
\]
Then we have
\[
\tr_p\circ S(x)=\psi(x)
\quad\text{for}\
x\in\bM_m.
\]
Moreover, $S$ is faithful.
By Lemma \ref{lem:3.1},
there exists a unique ucp map
$T\colon\bM_p\to\bM_m$ 
such that
\[
\tr_m(x^*h^{1/2}T(y)h^{1/2})=\tr_p(S(x)^*y)
\quad\text{for}\
x\in\bM_m, y\in\bM_p.
\]
In particular, we have
\[
\psi\circ T(y)
=\tr_p(y)
\quad\text{for}\
y\in\bM_p.
\]
Moreover, for $1\leq k,l\leq m$, we have
\[
\tr_m(e_{kl}^*h^{1/2}T(y)h^{1/2})
=\tr_p(S(e_{kl})^*y)
=\tr_p(f_{kl}^*y)
\quad\text{for}\
y\in\bM_p.
\]
Hence the $(k,l)$-th element
of the $m\times m$ matrix $h^{1/2}T(y)h^{1/2}$
is $m\tr_p(f_{kl}^*y)$.
This implies that
the $(k,l)$-th element
of the $m\times m$ matrix $T(y)$
is
\[
m\lambda_k^{-1/2}\lambda_l^{-1/2}\tr_p(f_{kl}^*y)
=
pp_k^{-1/2}p_l^{-1/2}\tr_p(f_{kl}^*y).
\]

Note that the $(i,j)$-th element 
of the $m\times m$ matrix
$T\circ S(e_{ij})=T(f_{ij})$
is
\[
(p_ip_j)^{-1/2}\min\{p_i, p_j\}
\]
and all other elements of the matrix are zero.
Hence
\[
T\circ S(e_{ij})=(p_ip_j)^{-1/2}\min\{p_i, p_j\}e_{ij}.
\]

Therefore we can obtain
\[
(\|x-T\circ S(x)\|_\psi^\sharp)^2
\leq
\|h^{1/2}x-xh^{1/2}\|_2^2.
\]
\epf

\brem\label{rem:commute}
Let $k$ be an $m\times m$-diagonal matrix
with eigenvalues
$\mu_1,\dots,\mu_m$.
Let $p$ be the positive integer and $T$ the ucp map as in the proof of Lemma \ref{lem:3.2}.
We define a $p\times p$-diagonal matrix
$\bar{k}$ by
\[
\bar{k}\coloneqq \sum_{i=1}^m\mu_iF_{ii}\in\bM_p
\]
Then we have
\[
kT(y)=T(\bar{k}y)
\quad\text{and}\quad
T(y)k=T(y\bar{k})
\quad\text{for}\
y\in \bM_p.
\]
\erem

The third lemma is also given in \cite{ha1}
in the case where $\varphi$ is tracial.

\blem[cf.\ {\cite[Lemma 3.3]{ha1}}]\label{lem:3.3}
Let $T\colon\bM_m\to M$
be a ucp map.
For $\e>0$,
there exists a ucp\ map 
$T'\colon\bM_m\to M$
such that 
$\|T-T'\|<\e$ and
\[
\varphi\circ T'(x)=\tr_m(xh')
\quad\text{for}\
x\in\bM_m,
\]
where 
$h'\in\bM_m^+$
has strictly positive rational eigenvalues.
\elem

Now we prove the main theorem in this section.

\bthm[cf.\ {\cite[Theorem 3.1]{ha4}}]\label{thm:semidiscrete}
If $M$ is injective,
then
for any $u_1,\dots,u_n\in\cU(M)$, any $\e>0$ and $\delta>0$,
there exist a ucp map
$T\colon\bM_m\to M$
and 
$v_1,\dots,v_n\in\cU(\bM_m)$
such that
$\psi=\varphi\circ T$ is a fn state 
on $\bM_m$, and
\begin{align*}
&\|\sigma_t^\varphi\circ T-T\circ\sigma_t^\psi\|
\leq\delta|t|
\quad\text{for}\
t\in\R,
\\
&\|T(v_k)-u_k\|_\varphi<\e
\quad\text{for}\
1\leq k\leq n.
\end{align*}
\ethm
\bpf


{\bf Step 0} : 


The first part is the same as in
\cite[Lemma 3.4]{ha4}.
Let $N\coloneqq M\rtimes_{\sigma^\varphi}\R$.
We denote by $\lambda^\varphi(t)$ 
the implementing unitary for $\sigma_t^\varphi$,
and by $\theta^\varphi$ the dual action of $\sigma^\varphi$.
Then there exists a fns (faithful normal semifinite) operator-valued weight
$P\colon N^+\to\widehat{M}^+$, which is given by
\[
P(y)\coloneqq \int_\R\theta_s^\varphi(y)ds
\quad\text{for}\
y\in N^+.
\]
Let $\widetilde{\varphi}\coloneqq\varphi\circ P$ be the dual weight of $\varphi$.
Recall that $N$ has a fns trace $\tau$
such that
\[
\tau\circ\theta_s^\varphi=e^{-s}\tau
\quad
\text{for}\
s\in\R.
\]
Let $a$ be the positive self-adjoint operator 
affiliated with $N_{\widetilde{\varphi}}$
such that $\exp(ita)=\lambda^\varphi(t)$ for $t\in\R$.
Then $\tau$ is given by
\[
\tau=\widetilde{\varphi}(e^{-a}\ \cdot\ ).
\]

Put $e_\alpha\coloneqq 1_{[0,\alpha]}(a)$ for $\alpha>0$.
Thanks to \cite[Lemma 3.4]{ha4}, 
we have 
$P(e_\alpha)=\alpha 1$.
Hence 
$\widetilde{\varphi}(e_\alpha)=\alpha$
and $\tau(e_\alpha)=1-e^{-\alpha}<\infty$.
By using \cite[Lemma 3.4]{ha4} again, we obtain
\[
\lim_{\alpha\to\infty}\left\|
\frac{1}{\alpha}P(e_\alpha x e_\alpha)-x
\right\|_\varphi=0
\quad\text{for}\
x\in M.
\]

Let $u_1,\dots,u_n\in\cU(M)$,
$\e>0$ and $\delta>0$ be given.
Take $\e'>0$ such that
\[
(2\e')^{1/2}+\e'<\e.
\]
Then we choose $\e_0,\e_1,\e_2, \e_3>0$
such that $1>\e_3$ and
\[
8\e_3^{1/2}+\e_2+\e_1+\e_0<\e'.
\]

Take $\alpha>0$ such that
\[
\left\|
\frac{1}{\alpha}P(e_\alpha u_k e_\alpha)-u_k
\right\|_\varphi<\e_0
\quad\text{for}\
1\leq k\leq n.
\]
We define the ucp map 
$T_0\coloneqq  \alpha^{-1}P|_{e_\alpha Ne_\alpha} 
\colon e_\alpha Ne_\alpha\to M$
and the fn state $\varphi_0$ on $e_\alpha Ne_\alpha$ by
\[
\varphi_0\coloneqq \varphi\circ T_0
=\frac{1}{\alpha}\varphi\circ P(e_\alpha\ \cdot\ e_\alpha)
=\frac{1}{\alpha}\widetilde{\varphi}(e_\alpha\ \cdot\ e_\alpha).
\]
Set $x_k\coloneqq e_\alpha u_ke_\alpha\in e_\alpha Ne_\alpha$ for $1\leq k\leq n$.
Then $\|x_k\|\leq 1$ and
\[
\|T_0(x_k)-u_k\|_\varphi<\e_0
\quad\text{for}\
1\leq k\leq n.
\]
Moreover, we have 
\[
\sigma_t^\varphi\circ T_0=T_0\circ\sigma_t^{\varphi_0}
\quad
\text{for}\ 
t\in\R.
\]

We define the fn tracial state 
$\tau_\alpha$ on $e_\alpha Ne_\alpha$
by
\[
\tau_\alpha\coloneqq \frac{1}{1-e^{-\alpha}}\widetilde{\varphi}(e^{-a}e_\alpha\ \cdot\ )
\]
Set 
\[
h_0\coloneqq \frac{d\varphi_0}{d\tau_\alpha}=\frac{1-e^{-\alpha}}{\alpha}e^ae_\alpha,
\]
and $\varphi_0=\tau_\alpha(h_0\ \cdot\ )$.
Note that 
\[
\Sp(h_0)\subset[c^{-1}, c]
\quad
\text{for some}\
c>1.
\]


{\bf Step 1} : 

The second part is similar as in
\cite[Theorem 3.1]{ha4}.
Take $\delta'>0$ with $\delta'<\e_1$ such that 
if positive elements $a, b$
with $\Sp(a), \Sp(b)\subset [c^{-1}, c]$ and $\|a-b\|<\delta'$, 
then $\|\log(a)-\log(b)\|<\delta/2$.

Choose $(2c)^{-1}>\delta_1>0$ such that 
$3c^2\delta_1<\delta'$.
Then take $\lambda\in\Q$ 
such that $0<\lambda<1$, 
$\lambda^{-1}-1<\delta_1$.
Then $1-\lambda=\lambda(\lambda^{-1}-1)<\delta_1$.
Set
\[
J\coloneqq
\max\{j\in\N\mid \lambda^{-j}\leq c\}.
\]
Since $\lambda^{-J}\leq c$ and $c<\lambda^{-(J+1)}$, 
we have
$c^{-1}\leq \lambda^J$ and $\lambda^{J+1}<c^{-1}$.
By using the spectral decomposition of $h_0$,
we can choose projections 
\[
e_{J+1},e_J,\dots,e_1,e_0,e_{-1},\dots,e_{-J}
\]
with $\sum_{-J\leq j\leq J+1}e_j=1$
such that 
\[
h_0'\coloneqq\sum_{-J\leq j\leq J+1}\lambda^je_j\leq h_0,
\quad\text{and}\quad
\|h_0-h_0'\|<c\delta_1.
\]
Put $C\coloneqq\tau_\alpha(h_0')\leq 1$ 
and $h_1\coloneqq C^{-1}h_0'$.
Since $1-C=\tau_\alpha(h_0-h_0')\leq c\delta_1$,
we have 
\[
1\leq C^{-1}\leq(1-c\delta_1)^{-1}.
\]
Then
\begin{align*}
\|h_0-h_1\|
&\leq
\|h_0-h_0'\|+\|h_0'-h_1\|
\\
&\leq
c\delta_1+(C^{-1}-1)\|h_0'\|
\\
&\leq(1+2c)c\delta_1<3c^2\delta_1<\delta'.
\end{align*}
Hence we have
\[
\|\log(h_0)-\log(h_1)\|<\frac{\delta}{2}.
\]
Put the fn state
$\varphi_1\coloneqq \tau_\alpha(h_1\ \cdot\ )$ on $e_\alpha Ne_\alpha$.
We define a cp map $T_0'\colon e_\alpha Ne_\alpha\to M$ by
\[
T_0'(x)\coloneqq T_0(b^{1/2}xb^{1/2})
\quad\text{for}\ x\in e_\alpha Ne_\alpha,
\]
where $b\coloneqq h_0^{-1}h_0'\leq 1$.
Then $T_0'(1)=T_0(b)\leq 1$ and
\[
\varphi\circ T_0'(x)
=
\varphi_0(b^{1/2}xb^{1/2})
=
\tau_\alpha(h_0b^{1/2}xb^{1/2})
=
\tau_\alpha(h_0'x)
\quad\text{for}\
x\in e_\alpha Ne_\alpha.
\]
Next we define a cp map $T_1\colon e_\alpha Ne_\alpha\to M$ by
\[
T_1(x)\coloneqq T_0'(x)
+
\frac{\tau_\alpha((h_1-h_0')x)}{\tau_\alpha(h_1-h_0')}(1-T_0'(1))
\quad\text{for}\ x\in e_\alpha Ne_\alpha.
\]
Then $T_1(1)=1$ and
\[
\varphi\circ T_1(x)=\tau_\alpha(h_1x)=\varphi_1(x)
\quad\text{for}\ x\in e_\alpha Ne_\alpha.
\]
Moreover,
\begin{align*}
\|T_0(x)-T_1(x)\|
&\leq
\|x-b^{1/2}xb^{1/2}\|
+\|1-b\|\|x\|
\\
&=
\frac{1}{2}\|(1+b^{1/2})x(1-b^{1/2})+(1-b^{1/2})x(1+b^{1/2})\|
+\|1-b\|\|x\|
\\
&\leq
(\|1+b^{1/2}\|\|1-b^{1/2}\|+\|1-b\|)\|x\|
\\
&\leq 3c^2\delta_1\|x\|.
\end{align*}
Hence
\[
\|T_0-T_1\|\leq 3c^2\delta_1<\delta'<\e_1.
\]
Therefore
\[
\|T_1(x_k)-u_k\|_\varphi
\leq
\|T_1(x_k)-T_0(x_k)\|_\varphi+\|T_0(x_k)-u_k\|_\varphi
<\e_1+\e_0
\quad\text{for}\
1\leq k\leq n.
\]
Since $\sigma_t^\varphi\circ T_1=T_1\circ \sigma_t^{\varphi_0}$ and 
\[
\|h_0^{it}-h_1^{it}\|\leq\|\log(h_0)-\log(h_1)\||t|\leq\frac{\delta}{2}|t|,
\]
we have 
\begin{align*}
\|\sigma_t^{\varphi}\circ T_1(x)-T_1\circ\sigma_t^{\varphi_1}(x)\|
&=
\|T_1(\sigma_t^{\varphi_0}(x)-\sigma_t^{\varphi_1}(x))\|
\\
&\leq
\|\sigma_t^{\varphi_0}(x)-\sigma_t^{\varphi_1}(x)\|
\\
&=\|h_0^{it}xh_0^{-it}-h_1^{it}xh_1^{-it}\|
\\
&\leq
\delta|t|\|x\|.
\end{align*}
Hence 
\[
\|\sigma_t^{\varphi}\circ T_1-T_1\circ\sigma_t^{\varphi_1}\|
\leq\delta|t|
\quad\text{for}\
t\in\R.
\]


{\bf Step 2} : 


The next part is similar as in
\cite[Lemma 5.3]{ha2}.
For $j\in\N$,
we define the linear map
\[
E_j\colon e_\alpha Ne_\alpha\to N_j\coloneqq \{x\in e_\alpha Ne_\alpha \mid 
\sigma_t^{\varphi_1}(x)=\lambda^{ijt}x, t\in\R\}
\] 
by
\[
E_j(x)\coloneqq \frac{1}{t_0}\int_0^{t_0}\lambda^{-ijt}\sigma_t^{\varphi_1}(x)dt
\quad
\text{for}\
x\in e_\alpha Ne_\alpha,
\]
where $t_0\coloneqq -2\pi/\log\lambda$.
For $q\in\N$, we define the ucp\ map $\gamma_q$ on
$e_\alpha Ne_\alpha$ by 
\[
\gamma_q(x)\coloneqq \sum_{j=-q+1}^{q-1}\left(1-\frac{|j|}{q}\right)E_j(x)
\quad\text{for}\
x\in e_\alpha Ne_\alpha.
\]
By \cite[Lemma 5.2]{ha2}, 
$\varphi_1\circ\gamma_q=\varphi_1$
and $\|\gamma_q(x)-x\|_{\varphi_1}\to 0$ ($q\to\infty$)
for $x\in e_\alpha Ne_\alpha$.
Choose $q\in\N$ such that
\[
\|\gamma_q(x_k)-x_k\|_{\varphi_1}<\e_2
\quad\text{for}\
1\leq k\leq n.
\]

Let $(e_{rs})$ be the matrix units for $\bM_q$.
We define the fn state $\psi_\lambda$
on $\bM_q$ by
\[
\psi_\lambda\coloneqq \tr_q(h_\lambda\ \cdot\ ),
\]
where 
\[
h_\lambda\coloneqq \sum_{r=1}^q\lambda_re_{rr}
\quad\text{and}\quad
\lambda_r\coloneqq q(\sum_{s=1}^q\lambda^s)^{-1}\lambda^r.
\]
Note that
\[
\sigma_t^{\psi_\lambda}(e_{rs})
=h_\lambda^{it}(e_{rs})h_\lambda^{-it}
=\lambda^{i(r-s)t}e_{rs}.
\]
In particular, 
$\sigma_{t_0}^{\psi_\lambda}=\id$.
Put the fn state 
$\chi\coloneqq \varphi_1\otimes\psi_\lambda$ 
on $e_\alpha Ne_\alpha\otimes\bM_q$.
Since 
\[
\sigma_t^\chi(x\otimes e_{rs})=\lambda^{i(r-s)t}\sigma_t^{\varphi_1}(x)\otimes e_{rs},
\]
the centralizer $N_\chi\coloneqq (e_\alpha N e_\alpha\otimes\bM_q)_\chi$ is given by
\[
N_\chi=\left\{\left.
\sum_{r,s=1}^qx_{rs}\otimes e_{rs}
\in e_\alpha N e_\alpha\otimes\bM_q\
\right|\ x_{rs}\in N_{s-r}\right\}.
\]
We define ucp\ maps 
$S_2\colon e_\alpha N e_\alpha\to N_\chi$ 
and 
$T_2\colon N_\chi\to e_\alpha N e_\alpha$ 
by
\[
S_2(x)\coloneqq \sum_{r,s=1}^qE_{s-r}(x)\otimes e_{rs},
\]
\[
T_2\left(\sum_{r,s=1}^qx_{rs}\otimes e_{rs}\right)
\coloneqq \frac{1}{q}\sum_{r,s=1}^qx_{rs}.
\]
Since
$T_2\circ S_2=\gamma_q$,
we have
\[
\|T_2\circ S_2(x_k)-x_k\|_{\varphi_1}<\e_2
\quad\text{for}\
1\leq k\leq n.
\]

Put the state $\varphi_2\coloneqq \varphi_1\circ T_2$
on $N_\chi$.
Then $\varphi_2\circ S_2=\varphi_1\circ\gamma_q=\varphi_1$.
For $x_{rs}\in N_{s-r}$,
we have
\[
\varphi_2\left(
\sum_{r,s=1}^qx_{rs}\otimes e_{rs}
\right)=
\frac{1}{q}\sum_{r=1}^q\varphi_1(x_{rr}).
\]
Hence
\[
\varphi_2(y)=(\varphi_1\otimes\tr_q)(y)
\quad\text{for}\quad
y\in N_\chi\subset e_\alpha Ne_\alpha\otimes\bM_q.
\]
Namely $\varphi_2$ is the restriction of 
$\varphi_1\otimes\tr_q$
on $N_\chi$.
Hence we have
\[
\sigma_t^{\varphi_2}(y)
=(\sigma_t^{\varphi_1}\otimes\id)(y)
\quad\text{for}\
y\in N_\chi.
\]
By definition, we have
\[
\sigma_t^{\varphi_2}\circ S_2
=S_2\circ\sigma_t^{\varphi_1}
\quad\text{and}\quad
\sigma_t^{\varphi_1}\circ T_2
=T_2\circ\sigma_t^{\varphi_2}.
\]

Let $\tau_\chi$ be the restriction of $\chi$ on $N_\chi$,
which is tracial.
Then
\[
\tau_\chi(y)=(\varphi_1\otimes\psi_\lambda)(y)
=(\varphi_1\otimes\tr_q)((1\otimes h_\lambda)y)
\quad\text{for}\
y\in N_\chi.
\]
Since $1\otimes h_\lambda\in N_\chi$, we have
\[
\frac{d\varphi_2}{d\tau_\chi}=1\otimes h_\lambda^{-1},
\]
and
\[
\sigma_t^{\varphi_2}(y)
=
\left(\frac{d\varphi_2}{d\tau_\chi}\right)^{it}
y
\left(\frac{d\varphi_2}{d\tau_\chi}\right)^{-it}
=
(1\otimes h_\lambda^{-it})y(1\otimes h_\lambda^{it})
\quad\text{for}\
y\in N_\chi, t\in\R.
\]


{\bf Step 3} : 

In this step, we use the semidiscreteness of $M$.
Since $\|S_2(x_k)\|\leq 1$ for $1\leq k\leq n$,
we have
\[
S_2(x_k)=
\frac{1}{2}(w_{k1}+w_{k2})+\frac{i}{2}(w_{k3}+w_{k4})
\]
for some unitaries $w_{kl}$.
Put
\[
w_{kl}=\sum_{r,s}w_{rs}^{(kl)}\otimes e_{rs}
\in N_\chi\subset e_\alpha Ne_\alpha\otimes\bM_q
\quad\text{for}
\
1\leq k\leq n, 1\leq l\leq 4.
\]
Set $c_\lambda\coloneqq\max\{1, \|1\otimes h_\lambda^{-1}\|\}\geq 1$.
Since $e_\alpha Ne_\alpha$ is semidiscrete,
we can take ucp maps 
$S_3\colon e_\alpha Ne_\alpha\to \bM_p$
and 
$T_3\colon \bM_p\to e_\alpha Ne_\alpha$
such that
\[
\|T_3\circ S_3(w_{rs}^{(kl)})-w_{rs}^{(kl)}\|_{\varphi_1}<\frac{\e_3}{c_\lambda\sqrt{q}}
\quad\text{for}\
1\leq k\leq n, 1\leq l\leq 4, 1\leq r,s\leq q.
\]
Here, by using Lemma \ref{lem:3.3},
we may also assume that $T_3$ satisfies 
\[
\varphi_1\circ T_3(x)=\tr_p(h_3'x)
\quad\text{for}\
x\in\bM_p,
\]
where $h_3'\in\bM_p^+$ has strictly positive rational eigenvalues.
Let $E_\chi$ be the $\chi$-invariant conditional expectation 
from $e_\alpha Ne_\alpha\otimes\bM_q$ 
onto $N_\chi$, which is given by
\[
E_\chi\left(\sum_{r,s=1}^qx_{rs}\otimes e_{rs}\right)
=
\sum_{r,s=1}^qE_{s-r}(x_{rs})\otimes e_{rs}.
\]
Set $S_3^{(q)}\coloneqq S_3\otimes\id_{\bM_q}$
and $T_3^{(q)}\coloneqq T_3\otimes\id_{\bM_q}$.
Then we have
\[
\|E_\chi\circ T_3^{(q)}\circ S_3^{(q)}(w_{kl})
-w_{kl}\|_{\tau_\chi}
<
\frac{\e_3}{c_\lambda}
\]
and 
\[
\|T_3^{(q)}\circ S_3^{(q)}(S_2(x_k))
-S_2(x_k)\|_{\varphi_1\otimes\tr_q}
<2\e_3.
\]


{\bf Step 4} : 


Let $h_3\coloneqq h_3'\otimes h_\lambda\in\bM_p\otimes\bM_q=\bM_{pq}$
be the diagonal matrix with strictly positive rational eigenvalues.
We define
\[
\varphi_3\coloneqq \tau_\chi\circ E_\chi\circ T_3^{(q)}.
\]
Then 
\[
\varphi_3=\tr_{pq}(h_3\ \cdot\ ).
\]
By Lemma \ref{lem:3.2},
there exist ucp\ maps 
$S_4\colon \bM_{pq}\to \bM_m$
and 
$T_4\colon \bM_m\to \bM_{pq}$
such that
$\varphi_3\circ T_4=\tr_m$,
$\tr_m\circ S_4=\varphi_3$,
and
\[
\|T_4\circ S_4(y)-y\|_{\varphi_3}^\sharp\leq
\|yh_3^{1/2}-h_3^{1/2}y\|_2
\quad\text{for}\quad
y\in \bM_{pq}.
\]
Set $y_{kl}\coloneqq S_3^{(q)}(w_{kl})\in\bM_{pq}$.
By Lemma \ref{lem:3.1}, we have
\begin{align*}
\tr_{pq}(y_{kl}^*h_3^{1/2}y_{kl}h_3^{1/2})
&\geq
\|E_\chi\circ T_3^{(q)}(y_{kl})\|_{\tau_\chi}^2
\\
&\geq
(\|w_{kl}\|_{\tau_\chi}-\|E_\chi\circ T_3^{(q)}\circ S_3^{(q)}(w_{kl})-w_{kl}\|_{\tau_\chi})^2
\\
&>
\left(1-\frac{\e_3}{c_\lambda}\right)^2>1-\frac{2\e_3}{c_\lambda},
\end{align*}
Hence
\[
\|h_3^{1/2}y_{kl}-y_{kl}h_3^{1/2}\|_2^2
\leq2-2\tr_{pq}(y_{kl}^*h_3^{1/2}y_{kl}h_3^{1/2})
<
\frac{4\e_3}{c_\lambda}.
\]
Therefore
\[
\|T_4\circ S_4(y_{kl})-y_{kl}\|_{\varphi_3}^\sharp
<
\frac{2\e_3^{1/2}}{\sqrt{c_\lambda}}
\quad\text{for}\
1\leq k\leq n, 1\leq l\leq 4.
\]
Set 
\[
y_k\coloneqq \frac{1}{2}(y_{k1}+y_{k2})
+\frac{i}{2}(y_{k3}+y_{k4})
=S_3^{(q)}\circ S_2(x_k)\in\bM_{pq}
\quad\text{for}\
1\leq k\leq n.
\]
Then
\[
\|T_4\circ S_4(y_k)-y_k\|_{\varphi_3}^\sharp
\leq
\frac{4\e_3^{1/2}}{\sqrt{c_\lambda}}
\quad\text{for}\
1\leq k\leq n.
\]


{\bf Step 5} : 

Set $k_\lambda\coloneqq 1\otimes h_\lambda^{-1}\in
\bM_p\otimes\bM_q=\bM_{pq}$.
The $m\times m$-diagonal matrix $\bar{k}_\lambda$ is 
defined in Remark \ref{rem:commute}.
We define a fn state $\psi$ on $\bM_m$ by 
\[
\psi\coloneqq\tr_m(\bar{k}_\lambda\, \cdot\,).
\]
Now we define a ucp map
$T\colon \bM_m\to M$ by
\[
T\coloneqq T_1\circ T_2\circ E_\chi\circ T_3^{(q)}\circ T_4.
\]
Then by Remark \ref{rem:commute}, 
for $z\in\bM_m$, we have
\begin{align*}
\varphi\circ T(z)
&=
\varphi\circ T_1\circ T_2\circ E_\chi\circ T_3^{(q)}\circ T_4(z)
\\
&=
\varphi_1\circ T_2\circ E_\chi\circ T_3^{(q)}\circ T_4(z)
\\
&=
\varphi_2\circ E_\chi\circ T_3^{(q)}\circ T_4(z)
\\
&=
\tau_\chi(1\otimes h_\lambda^{-1}(E_\chi\circ T_3^{(q)}\circ T_4)(z))
\\
&=
\tau_\chi\circ E_\chi\circ T_3^{(q)}(k_\lambda T_4(z)))
\\
&=
\tau_\chi\circ E_\chi\circ T_3^{(q)}\circ T_4(\bar{k}_\lambda z)
\\
&=
\varphi_3\circ T_4(\bar{k}_\lambda z)
\\
&=
\tr_m(\bar{k}_\lambda z)=\psi(z).
\end{align*}
Hence $\psi=\varphi\circ T$, and
\[
\sigma_t^\psi(z)=\bar{k}_\lambda^{it}z\bar{k}_\lambda^{-it}
\quad\text{for}\
z\in\bM_m.
\]
By Remark \ref{rem:commute},  for $z\in\bM_m$ we have
\begin{align*}
\sigma_t^{\varphi_2}\circ
E_\chi\circ T_3^{(q)}\circ T_4(z)
&=
1\otimes h_\lambda^{-it}
(E_\chi\circ T_3^{(q)}\circ T_4(z))1\otimes h_\lambda^{it}
\\
&=
E_\chi\circ T_3^{(q)}(k_\lambda^{it}T_4(z)k_\lambda^{-it})
\\
&=
E_\chi\circ T_3^{(q)}\circ T_4(\bar{k}_\lambda^{it}z\bar{k}_\lambda^{-it})
\\
&=
E_\chi\circ T_3^{(q)}\circ T_4\circ\sigma_t^\psi(z).
\end{align*}
Therefore we have
\[
\|\sigma_t^\varphi\circ T-T\circ\sigma_t^\psi\|
\leq\delta|t|
\quad\text{for}\ t\in\R.
\]
Put $z_k\coloneqq S_4(y_k)
=S_4\circ S_3^{(q)}\circ S_2(x_k)\in\bM_m$
for $1\leq k\leq n$.
By Kadison-Schwarz inequality, we have
\begin{align*}
\|E_\chi\circ T_3^{(q)}\circ T_4(z_k)-
E_\chi\circ T_3^{(q)}(y_k)\|_{\varphi_2}^2
&\leq
\|k_\lambda\|\|T_4(z_k)-y_k\|_{\varphi_3}^2
\\
&\leq
2c_\lambda\|T_4(z_k)-y_k\|_{\varphi_3}^{\sharp 2}
\\
&\leq
32\e_3.
\end{align*}
Hence
\begin{align*}
&\|E_\chi\circ T_3^{(q)}\circ T_4(z_k)-S_2(x_k)\|_{\varphi_2}
\\
&\leq
\|E_\chi\circ T_3^{(q)}\circ T_4(z_k)-
E_\chi\circ T_3^{(q)}(y_k)\|_{\varphi_2}
+
\|E_\chi\circ T_3^{(q)}(y_k)
-E_\chi(S_2(x_k))\|_{\varphi_2}
\\
&\leq
\sqrt{32}\e_3^{1/2}+
\|T_3^{(q)}\circ S_3^{(q)}(S_2(x_k))
-S_2(x_k)\|_{\varphi_1\otimes\tr_q}
\\
&\leq 6\e_3^{1/2}+2\e_3\leq 8\e_3^{1/2}.
\end{align*}
Moreover
\begin{align*}
\|
T_1\circ T_2(
E_\chi\circ T_3^{(q)}\circ T_4(z_k)
-S_2(x_k))
\|_\varphi^2
&\leq
\|E_\chi\circ T_3^{(q)}\circ T_4(z_k)
-S_2(x_k)\|_{\varphi_2}^2
\\
&<64\e_3.
\end{align*}
Similarly we have
\[
\|
T_1(T_2\circ S_2(x_k)-x_k) 
\|_\varphi^2
\leq
\|T_2\circ S_2(x_k)-x_k\|_{\varphi_1}^2
<\e_2^2.
\]
Therefore
\begin{align*}
\|T(z_k)-u_k\|_\varphi
&\leq
\|
T_1\circ T_2(
E_\chi\circ T_3^{(q)}\circ T_4(z_k)
-S_2(x_k))
\|_\varphi
\\
&+
\|
T_1(T_2\circ S_2(x_k)-x_k) 
\|_\varphi
+
\|
T_1(x_k)-u_k
\|_\varphi
\\
&<
8\e_3^{1/2}+\e_2+\e_1+\e_0
<\e'
\end{align*}
By the polar decomposition,
we obtain unitaries $v_k$ in $\bM_m$ such that
\[
z_k= v_k|z_k|
\quad\text{for}\
1\leq k\leq n.
\]
Then 
\[
\|z_k\|_\psi^2
\geq
\|T(z_k)\|_\varphi^2
>1-2\e'.
\]
Since
$\|z_k\|\leq 1$
and
$|z_k|^2+(1-|z_k|)^2\leq 1$,
we have
\begin{align*}
\|v_k-z_k\|_\psi^2
&=
\|1-|z_k|\|_\psi^2
\\
&\leq
1-\||z_k|\|_\psi^2
\\
&<
2\e'.
\end{align*}
Therefore
\begin{align*}
\|T(v_k)-u_k\|_\varphi
&\leq
\|T(v_k-z_k)\|_\varphi+\|T(z_k)-u_k\|_\varphi
\\
&\leq
\|v_k-z_k\|_\psi+\e'
\\
&<(2\e')^{1/2}+\e'<\e.
\end{align*}

\epf

\brem\label{rem:semidiscrete}
In Theorem \ref{thm:semidiscrete},
we obtain the ucp map
$T\colon\bM_m\to M$
with $\psi=\varphi\circ T$
such that $h_\psi$ has strictly positive rational eigenvalues,
where $\psi=\tr_m(h_\psi\ \cdot\ )$.

Moreover, in Theorem \ref{thm:semidiscrete}, 
we assume $M$ is a finite von Neumann algebra 
with a fn tracial state $\varphi=\tau$.
Then we can choose the ucp map $T$ satisfying $\tau\circ T=\tr_m$.
This fact is exactly \cite[Lemma 3.4]{ha1}.
Indeed, in the proof of Theorem \ref{thm:semidiscrete},
we omit Step 0,1 and 2.
In Step 3,
we set $M=e_\alpha Ne_\alpha=N_\chi$,
$\tau=\varphi_1$,
$u_k=S_2(x_k)$ and $q=1$.
Then we obtain ucp maps 
$S_2\colon M\to\bM_p$ and $T_2\colon \bM_p\to M$
 such that 
 \[
 \|T_2\circ S_2(u_k)-u_k\|_\tau<\e_3
 \quad\text{for}\
 1\leq k\leq n.
 \]
 In Step 4,
 we set $\varphi_3=\tau\circ T_2$
 and $y_k=S_2(u_k)$.
 Then we have ucp maps
$S_3\colon \bM_p\to\bM_m$ and $T_3\colon \bM_m\to \bM_p$
 such that 
\[
\|T_3\circ S_3(y_k)-y_k\|_{\varphi_3}<2\e_3^{1/2}
\quad\text{for}\
1\leq k\leq n.
\]
In Step 5, if we define $T\coloneqq T_2\circ T_3$
and $z_k=S_3(y_k)$,
then we have $\tau\circ T=\tr_m$ and 
\[
\|T(z_k)-u_k\|_\tau\leq\e_3+2\e_3^{1/2}
\quad\text{for}\ 1\leq k\leq n.
\] 

Next we consider the case where 
$M$ is a type III$_\lambda$ factor with $0<\lambda<1$
and a fn state $\varphi$ on $M$ satisfies 
$\sigma_{t_0}^\varphi=\id$
with $t_0=-2\pi/\log\lambda$.
Then we can choose the ucp map $T\colon\bM_m\to M$
with $\psi=\varphi\circ T=\tr_m(h_\psi\ \cdot\ )$
such that 
\[
\lambda_1/\lambda_2\in\{\lambda^n\mid n\in\Z\}
\quad\text{for}\
\lambda_1,\lambda_2\in\Sp(h_\psi).
\]
This fact is weaker than \cite[Theorem 3.4]{ha2},
but it is sufficient for our purpose.
By identifying $\widehat{\Z}$ with $\R/t_0\Z$,
\[
N_0=M\rtimes_{\sigma^\varphi}(\R/t_0\Z).
\]
is generated by
$\pi_0^\varphi(x)$ and $\lambda_0^\varphi(t)$,
where
\begin{align*}
(\pi_0^\varphi(x)\xi)(s)
&=
\sigma_{-s}^\varphi(x)\xi(s) \\
(\lambda_0^\varphi(t)\xi)(s)
&=
\xi(s-t)
\end{align*}
for $\xi\in L^2(\R/t_0\Z, H_\varphi)$.
By \cite[Proposition 5.6]{hs},
we have $N\simeq N_0\otimes L^\infty(0, \gamma_0)$
by identifying
\begin{align*}
\pi^\varphi(x)&=\pi_0^\varphi(x)\otimes 1 \\
\lambda^\varphi(t)&=\lambda_0^\varphi(t)\otimes m(e^{it}),
\end{align*}
where $\gamma_0=-\log\lambda$
and $m(e^{it})$ is the multiplication operator
\[
(m(e^{it})\xi)(\gamma)=e^{it\gamma}\xi(\gamma)
\quad\text{for}\
\xi\in L^2(0, \gamma_0).
\]
We denote the canonical traces $\tau_0$ and $\tau$,
the dual weights $\widetilde{\varphi}_0$ and $\widetilde{\varphi}$,
on $N_0$ and $N$, respectively.
Let $h_\varphi=d\widetilde{\varphi}_0/d\tau_0$
and $k_\varphi=d\widetilde{\varphi}/d\tau$.
Then
\[
k_\varphi=h_\varphi\otimes m(e^\gamma).
\]
Note that 
$\lambda_0^\varphi(t+t_0)=\lambda_0^\varphi(t)$
for $t\in\R$ and $h_\varphi^{it}=\lambda_0^\varphi(t)$.
Hence $\Sp(h_\varphi)=\{\lambda^n\}_{n\in\Z}\cup\{0\}$.
Therefore we have
\[
h_0=\frac{1-e^{-\alpha}}{\alpha}(h_\varphi\otimes m(e^\gamma))e_\alpha
\]
In Step 1,
if $\lambda\not\in\Q$, then
we can take $\mu\in\Q$ such that $\mu$ is sufficiently close to $\lambda$ and define
\[
h_1\coloneqq\frac{1-e^{-\alpha}}{\alpha}(h_\varphi'\otimes m(e^\gamma))e_\alpha.
\]
satisfying $\|h_0-h_1\|<\delta'$
and $\Sp(h_\varphi')\subset\{\mu^n\}_{n\in\Z}\cup\{0\}$.
Therefore by the proof of Theorem \ref{thm:semidiscrete},
we have the ucp map 
$T\colon\bM_m\to M$
with $\psi=\varphi\circ T$ 
such that 
$\psi=\tr_m(\bar{k}_\mu\ \cdot\ )$.
Then by small perturbation of $T$, 
we can obtain $T'$ such that 
$\psi=\varphi\circ T'=\tr_m(\bar{k}_\lambda\ \cdot\ )$ .
\erem


\section{The bicentralizer of a type III$_0$ factor}


Let $M$ be a von Neumann algebra.
We denote by $\cW(M)$ and $\cW_0(M)$
the set of ns (normal semifinite) weights 
and fns (faithful normal semifinite) weights on $M$,
respectively.
For $\varphi\in\cW(M)$, we define 
\[
\mathfrak{n}_\varphi\coloneqq\{x\in M \mid \varphi(x^*x)<\infty\}
\]
and $\mathfrak{m}_\varphi\coloneqq\mathfrak{n}_\varphi^*\mathfrak{n}_\varphi$.
Let $\omega\in\beta(\N)\setminus\N$
be a free ultrafilter on $\N$.

\bdf[{\cite[Definition 3.1]{hi}}]
For a fn state $\varphi$ on $M$,
we define 
the {\em asymptotic centralizer}
(resp.
{\em $\omega$-asymptotic centralizer}) 
of $\varphi$ by
\[
\mathrm{AC}\, (M,\varphi)
\coloneqq \left\{
(x_n)_n\in\ell^\infty(\N, M)
\mid
\lim_{n\to\infty}\|x_n\varphi-\varphi x_n\|=0
\right\}.
\]
\[
\mathrm{AC}_\omega (M,\varphi)
\coloneqq \left\{
(x_n)_n\in\ell^\infty(\N, M)
\mid
\lim_{n\to\omega}\|x_n\varphi-\varphi x_n\|=0
\right\}.
\]
We also define 
the {\em bicentralizer}
(resp.
{\em $\omega$-bicentralizer}) 
of $\varphi$ by
\[
\mathrm{B}\,(M, \varphi)
\coloneqq \left\{
a\in M
\mid
\lim_{n\to\infty}\|ax_n-x_na\|_\varphi=0
\quad\text{for}\
(x_n)_n\in\mathrm{AC}_\omega (M,\varphi)
\right\}.
\]
\[
\mathrm{B}_\omega(M, \varphi)
\coloneqq \left\{
a\in M
\mid
\lim_{n\to\omega}\|ax_n-x_na\|_\varphi=0
\quad\text{for}\
(x_n)_n\in\mathrm{AC}_\omega (M,\varphi)
\right\}.
\]
\edf
 
 We define
 \begin{align*}
 \mathcal{I}_\omega(M)
 &\coloneqq \left\{
 (x_n)_n\in\ell^\infty(M)
 \mid
 x_n\to0
 \ \text{$*$-strongly as}\ n\to\omega
 \right\}\\
 \mathcal{M}^\omega(M)
 &\coloneqq \left\{
 (x_n)_n\in\ell^\infty(M)
 \mid
 (x_n)_n\mathcal{I}_\omega(M)\subset\mathcal{I}_\omega(M)
\ \text{and}\ 
 \mathcal{I}_\omega(M) (x_n)_n\subset\mathcal{I}_\omega(M)
 \right\}
 \end{align*}
Then the {\em multiplier algebra} $\mathcal{M}^\omega(M)$
is a C$^*$-algebra and 
$\mathcal{I}_\omega(M)\subset\mathcal{M}^\omega(M)$
is a norm closed two-sided ideal.
Following \cite{oc},
we define the {\em ultrapower von Neumann algebra} $M^\omega\coloneqq \mathcal{M}^\omega(M)/\mathcal{I}_\omega(M)$, which is indeed well-known to be a von Neumann algebra.

\bdf[{\cite[Definition 4.25]{ah}}]
For $\varphi\in\cW(M)$,
we define $\varphi^\omega\in\cW(M)$
by
\[
\varphi^\omega\coloneqq \varphi\circ E,
\]
where 
\[
E\colon M^\omega\ni(x_n)^\omega
\mapsto{\text{wot-}}\lim_{n\to\omega}x_n\in M
\]
is the canonical fn conditional expectation.
If $\varphi$ is faithful, then so is $\varphi^\omega$.
\edf

The following fact induces us to define the bicentralizer of a general weight.

\bprop[{\cite[Proposition 3.2, Proposition 3.3]{hi}}]\label{prop:hi}
For a fn state $\varphi$ on $M$, we have
\[
\mathrm{B}\,(M, \varphi)
=\mathrm{B}_\omega(M,\varphi)
=[(M^\omega)_{\varphi^\omega}]'\cap M.
\]
In particular, the bicentralizer of $\varphi$
does not depend on the choice of a free ultrafilter 
$\omega\in\beta(\N)\setminus \N$.
\eprop

\bdf
For any $\varphi\in\cW_0(M)$,
we define the $\omega$-{\em bicentralizer} 
of $\varphi$ by
\[
B_\omega(M, \varphi)
\coloneqq [(M^\omega)_{\varphi^\omega}]'\cap M.
\]
\edf

\brem\label{quest:omega}
Let $\varphi\in\cW_0(M)$.
Does $B_\omega(M, \varphi)$ depend on the choice
of a free ultrafilter
$\omega\in\beta(\N)\setminus \N$?
We give a partial answer of the question
in the end of this section.
\erem

\blem\label{lem:lacunary}
If $\varphi\in\cW_0(M)$ is lacunary,
then we have
\[
B_\omega(M, \varphi)\supset\cZ(M_\varphi).
\]
\elem

\bpf
Since 
$
(M^\omega)_{\varphi^\omega}=(M_\varphi)^\omega
$
by \cite[Proposition 4.27]{ah}, 
we have
\[
B_\omega(M, \varphi)
=
[(M^\omega)_{\varphi^\omega}]'\cap M
=
[(M_\varphi)^\omega]'\cap M
\supset\cZ(M_\varphi).
\]
\epf

\brem
If $M$ is a type III$_\lambda$ factor
with $0\leq \lambda<1$
and $\varphi\in\cW_0(M)$ is lacunary, then
by \cite[Lemma XII.4.7]{tak},
we have
\[
B_\omega(M, \varphi)
=
[(M^\omega)_{\varphi^\omega}]'\cap M
=
[(M_\varphi)^\omega]'\cap M
\subset M_\varphi'\cap M=\cZ(M_\varphi),
\]
and therefore
\[
B_\omega(M, \varphi)=\cZ(M_\varphi).
\]
\erem

The following arguments are based on the work of Connes-Takesaki \cite{ct}.
From now on, we assume that $M$
is a $\sigma$-finite type III$_\lambda$ factor,
for $0\leq\lambda<1$
and $\varphi\in\cW_0(M)$
is a lacunary weight of infinite multiplicity.
Then $M_\varphi$ is a type II$_\infty$ von Neumann algebra,
and there exists a unitary $U\in M$
such that 
\[
UM_\varphi U^*=M_\varphi,
\qquad
\varphi\circ\Ad(U)\leq\lambda_0\varphi,
\quad 0<\lambda_0<1;
\]
\[
M\simeq M_\varphi\rtimes_\theta\Z
\qquad
\theta=\Ad(U)|_{M_\varphi}.
\]
Moreover $\tau\coloneqq \varphi|_{M_\varphi}$
is a fns trace on $M_\varphi$.
For $m\in \Z$,
there exists non-singular positive self-adjoint operator
$\rho_m$ affiliated to $\cZ(M_\varphi)$
such that 
\benu
\item
$\varphi\circ\theta^m=\varphi_{\rho_m}$,
\item
$\rho_{m+n}=\rho_m\theta^{-m}(\rho_n)$,
\item
$\sigma_t^\varphi(U^m)=U^m\rho_m^{it}$,
\item
$\rho_m\leq\rho_1<1$ for $m>0$,
\item
$\rho_m\geq\rho_{-1}>1$ for $m<0$,
\eenu
We simply write $\rho\coloneqq \rho_1$.

\brem\label{rem:inter}
Let $u, v\in M$
be unitaries and $x\in M$ with the polar decomposition $x=w|x|$.
If $ux=xv$,
then $uw=wv$,
$u(ww^*)=(ww^*)u$
and $v(w^*w)=(w^*w)v$.
\erem

\blem[cf. {\cite[\S 3.2, Lemma 2.6]{ct}}]\label{lem:ct}
Suppose that $\psi_j=\varphi_{h_j}$
with $h_j\in M_\varphi^+$ satisfying 
\[
\rho s(h_j)\leq h_j<1
\quad \text{for} \quad j=1, 2.
\]
It $\psi_2^\omega=(\psi_1^\omega)_u$
for a partial isometry $u\in M^\omega$
with $uu^*=s(\psi_1^\omega)$ 
and $u^*u=s(\psi_2^\omega)$,
then $u\in (M^\omega)_{\varphi^\omega}=(M_\varphi)^\omega$.
\elem

\bpf
The proof is similar to the one of \cite[Lemma XII.4.14]{tak}.
Set 
\[
k_j\coloneqq \rho(1-s(h_j))+h_j
\quad
\text{for}
\
j=1, 2.
\]
Then 
$\rho\leq k_j<1$.
Note that 
\[
s(\psi_j^\omega)=s(\psi_j)=s(h_j)\in M_\varphi
\]
by \cite[Lemma XII.4.13]{tak}.

By {\cite[Lemma XII.4,3]{tak}},
we have
\begin{align*}
uk_2^{it}
&=
us(h_2)k_2^{it}
=
uh_2^{it}
=
u(D\psi_2^\omega \colon D\varphi^\omega)_t
=
uu^*(D\psi_1^\omega \colon D\varphi^\omega)_t\sigma_t^{\varphi^\omega}(u)
\\
&=
s(h_1)h_1^{it}\sigma_t^{\varphi^\omega}(u)
=
k_1^{it}s(h_1)\sigma_t^{\varphi^\omega}(u)
=
k_1^{it}\sigma_t^{\varphi^\omega}(u).
\end{align*}
Therefore we obtain 
\[
uk_2^{it}=k_1^{it}\sigma_t^{\varphi^\omega}(u)
\quad 
\text{for}
\quad
t\in \R.
\]
By \cite[Proposition 6.23]{ah},
$M^\omega$ is canonically isomorphic to 
$(M_\varphi)^\omega\rtimes_{\theta^\omega}\Z$.
Hence we choose a sequence$(x^{(m)})$ in $(M_\varphi)^\omega=(M^\omega)_{\varphi^\omega}$.
such that $u=\sum_{m\in\Z}x^{(m)}U^m$ in $M^\omega$.
Now we have
\[
uk_2^{it}
=
\sum_{m\in\Z}x^{(m)}U^mk_2^{it}
=
\sum_{m\in\Z}x^{(m)}\theta^m(k_2^{it})U^m,
\]
and
\[
k_1^{it}\sigma_t^{\varphi^\omega}(u)
=
k_1^{it}\sum_{m\in\Z}\sigma_t^{\varphi^\omega}(x^{(m)}U^m)
=
\sum_{m\in\Z}k_1^{it}x^{(m)}U^m\rho_m^{it}
=
\sum_{m\in\Z}k_1^{it}x^{(m)}\theta^m(\rho_m^{it})U^m.
\]
By the uniqueness of the expansion, we have
\[
k_1^{it}x^{(m)}\theta^m(\rho_m^{it})
=
x^{(m)}\theta^m(k_2^{it}).
\]
Hence
\[
k_1^{it}x^{(m)}
=
x^{(m)}\theta^m(k_2^{it}\rho_m^{-it}).
\]
For each $m\in\Z$, by Remark \ref{rem:inter},
we may and do assume that 
$w\coloneqq x^{(m)}$ is a partial isometry in $(M_\varphi)^\omega$
such that
$w^*w$ commutes with $\theta^m(k_2^{it}\rho_m^{-it})$
and $ww^*$ commutes with $k_1^{it}$.
Then
\[
k_1^{it}ww^*=w\theta^m(k_2^{it}\rho_m^{-it})w^*.
\]

If $m>0$, then $\rho_m\leq\rho\leq k_2$.
Hence $H\coloneqq \theta^m(k_2\rho_m^{-1})\geq 1$ and $0\leq K\coloneqq k_1< 1$.
Then the functions
\[
\{z\in \C \mid \image(z)\geq 0\}\ni z\mapsto wH^{iz}w^*
\]
and
\[
\{z\in \C \mid \image(z)\leq 0\}\ni z\mapsto K^{iz}ww^*
\]
are analytic and bounded,
and $H^{it}ww^*=wK^{it}w^*$ for $t\in\R$.
Hence the function $z\mapsto K^{iz}ww^*$ 
can be extended to a bounded entire function,
which must be constant by Liouville's theorem.
Therefore $ww^*=ww^*k_1$ and $ww^*(1-k_1)=0$
Since $1-k_1$ is non-singular, we have $ww^*=0$.
Similarly we can show the case $m<0$.
\epf

\blem[cf. {\cite[Lemma XII.4.14]{tak}}]\label{lem:central}
Let $\psi\in\cW_0(M)$.
If $\psi=\varphi_h$ for some
$h\in M_\varphi^+$ with
$\rho s(h)\leq h<1$,
then
\[
(M^\omega)_{\psi^\omega}
\subset (M^\omega)_{\varphi^\omega}
=(M_\varphi)^\omega.
\]
\elem

\bpf
If $u\in (M^\omega)_{\psi^\omega}$ is a partial isometry
such that $uu^*=s(\psi^\omega)=u^*u$,
then $\psi^\omega=\psi_u^\omega$.
By Lemma \ref{lem:ct}, we have 
$u\in (M^\omega)_{\varphi^\omega}
=(M_\varphi)^\omega$.
\epf

\bthm\label{thm:bicentral}
If $M$ is a type III$_0$ factor
with separable prequel,
then 
\[
\mathrm{B}_\omega(M, \psi)\ne\C 1
\]
for any $\psi\in\cW_0(M)$. 
In particular, 
\[
\mathrm{B}\,(M, \psi)\ne\C 1
\]
for any fn state $\psi$ on $M$.
\ethm

\bpf
By \cite[Theorem XII.4.10]{tak},
there exists $h\in M_\varphi^+$
satisfying $\rho s(h)\leq h<1$
such that $\psi\sim\varphi_h$, i.e., 
there exists an isometry $u\in M$
such that $1=s(\psi)=u^*u$, $s(\varphi_h)=uu^*$ 
and $\psi(x)=\varphi_h(uxu^*)$ for $x\in M$.

By Lemma \ref{lem:central}, we have
\[
(M^\omega)_{\varphi_h^\omega}
\subset (M^\omega)_{\varphi^\omega}
=(M_\varphi)^\omega.
\]
If $a\in M_{\psi^\omega}$ satisfies $uau^*=0$,
then
\[
0=\varphi_h^\omega(ua^*u^*uau^*)
=\psi^\omega(a^*a).
\]
Since $\psi^\omega$ is faithful,
we have $a=0$.
Hence the adjoint map 
$\Ad(u)\colon M_{\psi^\omega}\to M_{\varphi_h^\omega}$
is an injective normal $*$-homomorphism.
Since
\[
u(M^\omega)_{\psi^\omega}u^*
\subset
(M^\omega)_{\varphi_h^\omega}
\subset
(M^\omega)_{\varphi^\omega}
\]
we have
\[
B_\omega(M, \varphi)
=
[(M^\omega)_{\varphi^\omega}]'\cap M
\subset
[(M^\omega)_{\varphi_h^\omega}]'\cap M
\subset
[u(M^\omega)_{\psi^\omega}u^*]'\cap M.
\]
By Lemma \ref{lem:lacunary},
we have 
\[
\cZ(M_\varphi)\subset B_\omega(M, \varphi)\subset [u(M^\omega)_{\psi^\omega}u^*]'\cap M.
\]
Let $a\in \cZ(M_\varphi)$.
For $x\in(M^\omega)_{\psi^\omega}$,
since $u$ is an isometry,
we have
\[
u^*aux=u^*a(uxu^*)u=u^*(uxu^*)au=xu^*au.
\]
Hence $u^*au\in B_\omega(M, \psi)$.

Now suppose that $u^*au\in\C 1$ 
for any $a\in\cZ(M_\varphi)$,
i.e., $u^*au=\gamma 1$ for some $\gamma\in\C$.
Recall that $\cZ(M_\varphi)=M_\varphi'\cap M$
by \cite[Corollary XII.4.17]{tak}.
Then $ae=\gamma e$,
because $e\coloneqq uu^*=s(h)\in M_\varphi$.
Since $axe=xae=\gamma xe$ for any $x\in M_\varphi$,
we have $ac(e)=\gamma c(e)$,
where $c(e)$ is the central support of $e$ in $M_\varphi$.
Therefore we have 
$\cZ(M_\varphi)c(e)=\C c(e)$,
which contradicts the fact that 
$\cZ(M_\varphi)$ is non-atomic
by \cite[Corollary XII.3.15]{tak}.
Hence it follows that $B_\omega(M, \psi)\ne\C 1$.
\epf

Next we discuss a von Neumann algebra $M$
with trivial bicentralizer, except type III$_1$.

\bprop\label{prop:trivialBC}
Let $\psi$ be a fn state on $M$
with $B(M, \psi)=\C 1$.
\benu
\item If $M$ is a semifinite von Neumann algebra,
then $M$ is finite and $\psi$ is tracial.
\item If $M$ is a type III$_\lambda$ factor with $0<\lambda<1$, 
then $\sigma_{t_0}^\psi=\id$, where $t_0=-2\pi/\log\lambda$.
\eenu
\eprop

\bpf
(1) Assume $M$ is a semifinite von Neumann algebra
with a fsn trace $\tau$.
Then there exists a non-singular positive self-adjoint operator $h$
affiliated with $M_\psi$
such that $\psi=\tau_h$
by \cite[Theorem VIII.3.14]{tak}.
Hence
\[
\sigma_t^\psi(x)=h^{it}xh^{-it}
\quad\text{for}\
t\in\R, x\in M.
\]
For $(x_n)^\omega\in(M^\omega)_{\psi^\omega}$,
thanks to \cite[Theorem 4.1]{ah}, we have
\[
h^{it}(x_n)^\omega h^{-it}
=
(\sigma_t^\psi(x_n))^\omega
=
\sigma_t^{\psi^\omega}((x_n)^\omega)
=(x_n)^\omega.
\]
Hence
\[
h^{it}\in [(M^\omega)_{\psi^\omega}]'\cap M=B(M,\psi)=\C1.
\]
Therefore we have $\sigma_t^\psi=\id$, which means that $\psi$ is tracial.

(2) We assume that $M$ is a type III$_\lambda$ factor 
with $0<\lambda<1$.
Let $\varphi\in\mathcal{W}_0(M)$ be 
a lacunary weight with infinite multiplicity.
By \cite[Corollary XII.4.10]{tak},
there exist $h\in M_\varphi^+$
and an isometry $u\in M$
such that
$1=s(\psi)=u^*u$, $s(\varphi_h)=uu^*$
and $\psi(x)=\varphi_h(uxu^*)$
for $x\in M$.
By Lemma \ref{lem:central}, we have
\[
(M^\omega)_{\varphi_h^\omega}
\subset(M^\omega)_{\varphi^\omega}=(M_\varphi)^\omega.
\]
By the proof of Theorem \ref{thm:bicentral},
\[
u(M^\omega)_{\psi^\omega}u^*
\subset (M^\omega)_{\varphi_h^\omega}.
\]
Hence we obtain
\[
h\in[(M^\omega)_{\varphi_h^\omega}]'\cap M
\subset
[u(M^\omega)_{\psi^\omega}u^*]'\cap M
\]
Then
\[
u^*hu\in[(M^\omega)_{\psi^\omega}]'\cap M
=B_\omega(M,\psi)=\C 1.
\]
Therefore $u^*hu=\gamma 1$
for some constant $\gamma$.
Since $h=uu^*huu^*=\gamma uu*=\gamma s(h)$,
for $x\in M$, we have
\[
\psi(x)=\varphi_h(uxu^*)
=\varphi(huxu^*)
=\varphi(uu^*huxu^*)
=\gamma\varphi(uxu^*).
\]
By \cite[Lemma XII.4.3]{tak},
\[
(D\psi\colon D\varphi)_t
=
u^*(D\gamma\varphi\colon D\varphi)_t\sigma_t^\varphi(u)
=
\gamma^{it}u^*\sigma_t^\varphi(u)
\eqqcolon u_t.
\]
Then
\[
\sigma_t^\psi(x)
=
u_t\sigma_t^\varphi(x)u_t^*
=
u^*\sigma_t^\varphi(u)\sigma_t^\varphi(x)\sigma_t^\varphi(u^*)u
=
u^*\sigma_t^\varphi(uxu^*)u.
\]
Therefore
\[
\sigma_{t_0}^\psi(x)
=
u^*\sigma_{t_0}^\varphi(uxu^*)u
=
u^*(uxu^*)u
=
x.
\]
\epf

Finally we discuss the problem in Remark \ref{quest:omega}.
Recall that $\varphi\in\cW_0(M)$ is {\em strictly semifinite}
if its restriction to $M_\varphi$ is also semifinite.

\bprop
If $\varphi\in\cW_0(M)$ is strictly semifinite,
then $B_\omega(M, \varphi)$ does not depend on the choice of a free ultrafilter $\omega\in\beta\N\setminus \N$.
\eprop

\bpf
We claim that 
$a\in B_\omega(M, \varphi)=[(M^\omega)_{\varphi^\omega}]'\cap M$ 
if and only if
\[
a\in\bigcap_{\substack{e\in M_\varphi \\ \varphi(e)<\infty}} 
[(M_e^\omega)_{\varphi_e^\omega}]'\cap M.
\]
Assume that $a\in M$ commutes with any element
in $(M_e^\omega)_{\varphi_e^\omega}$ for any projection $e\in M_\varphi$
with $\varphi(e)<\infty$.
Since $\varphi$ is strictly semifinite,
there exists an orthogonal family $(e_k)_{k\in I}$
of projections in $M_\varphi$ 
with sum $1$ such that $\varphi(e_k)<\infty$.
Put $p_F\coloneqq \sum_{k\in F}e_k$ for a finite subset $F\Subset I$.
Then $M_\varphi\ni p_F\nearrow 1$ and $\varphi(p_F)<\infty$.
Let $x\in (M^\omega)_{\varphi^\omega}$.
Then 
$p_Fxp_F\in(M_{p_F}^\omega)_{\varphi_{p_F}^\omega}$.
Hence 
\[
p_F(ax)p_F=a(p_Fxp_F)=(p_Fxp_F)a=p_F(xa)p_F
\]
and so $ax=xa$, namely $a\in B_\omega(M,\varphi)$.

Conversely, let $a\in B_\omega(M,\varphi)$.
Let $e\in M_\varphi$ be a projection with $\varphi(e)<\infty$.
Since $(M_e^\omega)_{\varphi_e^\omega}\subset (M^\omega)_{\varphi^\omega}$,
we have
$ax=xa$
for $x\in (M_e^\omega)_{\varphi_e^\omega}$.

Moreover,
\[
[(M_e^\omega)_{\varphi_e^\omega}]'\cap M=([(M_e^\omega)_{\varphi_e^\omega}]'\cap M_e)\oplus M_{e^\perp}.
\]
Indeed, let $a\in [(M_e^\omega)_{\varphi_e^\omega}]'\cap M$.
Since $e\in (M_e^\omega)_{\varphi_e^\omega}$,
we have 
\[
a=eae+eae^\perp\in ([(M_e^\omega)_{\varphi_e^\omega}]'\cap M_e)\oplus M_{e^\perp}.
\]
Conversely, Let $a=a_1+a_2\in ([(M_e^\omega)_{\varphi_e^\omega}]'\cap M_e)\oplus M_{e^\perp}$.
For any $x\in (M_e^\omega)_{\varphi_e^\omega}$,
we have $xa_2=a_2x=0$.
Hence $ax=xa$, and so $a\in [(M_e^\omega)_{\varphi_e^\omega}]'\cap M$.

By Proposition \ref{prop:hi},
we have 
\[
[(M_e^\omega)_{\varphi_e^\omega}]'\cap M_e=B_\omega(M_e, \varphi_e^\omega)
=B(M_e, \varphi_e^\omega).
\]
Therefore
\[
B_\omega(M, \varphi)=\bigcap_{\substack{e\in M_\varphi \\ \varphi(e)<\infty}} B(M_e, \varphi_e^\omega)\oplus M_{e^\perp},
\]
which means the independence on the choice of a free ultrafilter $\omega$.
\epf



\section{Almost unitary equivalence}


In this section, we generalize the notion of 
$\delta$-relatedness for two $n$-tuples of unit vectors
in a Hilbert bimodule.
In \cite[Remark 2.9]{ha4}, it is stated, 
but there is no proof, and so we enter into details.

Throughout this section, $M$ is a von Neumann algebra,
and $H$ is a Hilbert $M$-bimodule, i.e.,
$H$ is a Hilbert space with left and right actions 
\[
(x,\xi)\mapsto x\xi,
\quad
(x, \xi)\mapsto \xi x
\]
such that the above maps are bilinear and 
$(x\xi)y=x(\xi y)$
for $x, y\in M$, $\xi\in H$.
Moreover
\[
x\mapsto L_x,
\quad
x\mapsto R_x
\]
are a normal $*$-homomorphism and 
$*$-antihomomorphism, respectively,
where 
$L_x\xi\coloneqq x\xi$ 
and $R_x\xi\coloneqq\xi x$
for $x\in M$, $\xi\in H$.

\bdf[cf.\ {\cite[Definition 2.1]{ha2}}]\label{def:almost}
Two $n$-tuples 
$(\xi_1,\dots,\xi_n)$ and $(\eta_1,\dots,\eta_n)$
of unit vectors in $H$
are called {\em almost $\delta$-related} if 
for any $\e>0$,
there exist $a_1,\dots,a_p\in M$ such that
for $1\leq k\leq n$,
\begin{enumerate}
\item[(a1)] ${\dis \|(1-\sum_{j=1}^pa_j^*a_j)\xi_k\|<\e,
\|(1-\sum_{j=1}^pa_j^*a_j)\eta_k\|<\e}$;
\item[(a2)] ${\dis \|\xi_k(1-\sum_{j=1}^pa_j^*a_j)\|<\e,
\|\eta_k(1-\sum_{j=1}^pa_j^*a_j)\|<\e}$;
\item[(b1)] ${\dis \|(1-\sum_{j=1}^pa_ja_j^*)\xi_k\|<\e,
\|(1-\sum_{j=1}^pa_ja_j^*)\eta_k\|<\e}$;
\item[(b2)] ${\dis \|\xi_k(1-\sum_{j=1}^pa_ja_j^*)\|<\e,
\|\eta_k(1-\sum_{j=1}^pa_ja_j^*)\|<\e}$;
\item[(c)] ${\dis \sum_{j=1}^p\|a_j\xi_k-\eta_k a_j\|^2<\delta}$.
\end{enumerate}
\edf

\brem\label{rem:almost}
In this case, we can easily check that
\begin{itemize}
\item[(d)] ${\dis \sum_{j=1}^p\|a_j^*\eta_k-\xi_k a_j^*\|^2<2\delta}$ for $1\leq k\leq n$.
\end{itemize}
Indeed, 
for $\delta/4>\e>0$, we take $a_1,\dots,a_p\in M$ 
satisfying (a1),(a2),(b1),(b2) and (c).
Then we have
\begin{align*}
\sum_{j=1}^p\|a_j^*\eta_k-\xi_k a_j^*\|^2
&=
\sum_{j=1}^p\|a_j^*\eta_k\|^2+\|\xi_k a_j^*\|^2-2\real\ip{a_j^*\eta_k, \xi_k a_j^*} \\
&=
\ip{(\sum_{j=1}^pa_ja_j^*-1)\eta_k, \eta_k}
+\ip{\eta_k, \eta_k (1-\sum_{j=1}^pa_ja_j^*)}
\\
&\quad+\ip{\eta_k, \eta_k\sum_{j=1}^pa_ja_j^*} 
+\ip{\xi_k, \xi_k(\sum_{j=1}^pa_j^*a_j-1)}
\\
&\quad+\ip{(1-\sum_{j=1}^pa_j^*a_j)\xi_k, \xi_k}
+\ip{\sum_{j=1}^pa_j^*a_j\xi_k, \xi_k}
\\
&\quad-
2\real\ip{a_j\xi_k, \eta_k a_j} \\
&<4\e+\sum_{j=1}^p\|a_j\xi_k-\eta_k a_j\|^2
<2\delta.
\end{align*}
\erem

\brem\label{rem:almost2}
If two $n$-tuples 
$(\xi_1,\dots,\xi_n)$ and $(\eta_1,\dots,\eta_n)$
are almost $\delta$-related,
then for each $\e>0$,
we can choose $a_1,\dots,a_p$ such that for $1\leq k\leq n$,
\begin{enumerate}
\item[(a'1)] ${\dis \|(1-\sum_{j=1}^pa_j^*a_j)\xi_k\|<2\e,
\|(1-\sum_{j=1}^pa_j^*a_j)\eta_k\|<2\e}$;
\item[(a'2)] ${\dis \|\xi_k(1-\sum_{j=1}^pa_j^*a_j)\|<2\e,
\|\eta_k(1-\sum_{j=1}^pa_j^*a_j)\|<2\e}$;
\item[(b'1)] ${\dis \|(1-\sum_{j=1}^pa_ja_j^*)\xi_k\|<2\e,
\|(1-\sum_{j=1}^pa_ja_j^*)\eta_k\|<2\e}$;
\item[(b'2)] ${\dis \|\xi_k(1-\sum_{j=1}^pa_ja_j^*)\|<2\e,
\|\eta_k(1-\sum_{j=1}^pa_ja_j^*)\|<2\e}$;
\item[(c')] ${\dis \sum_{j=1}^p\|a_j\xi_k-\eta_k a_j\|^2<\delta}$;
\item[(d')] ${\dis \sum_{j=1}^p\|a_j^*\eta_k-\xi_k a_j^*\|^2<2\delta}$;
\item[(e')] ${\dis \|(\sum_{j=1}^pa_j^*a_j)\xi_k\|\leq 1,
\|(\sum_{j=1}^pa_j^*a_j)\eta_k\|\leq 1}$;
\item[(f')] ${\dis \|(\sum_{j=1}^pa_ja_j^*)\xi_k\|\leq 1,
\|(\sum_{j=1}^pa_ja_j^*)\eta_k\|\leq 1}$.
\end{enumerate}
Indeed, we take $a_1,\dots,a_p$ satisfying 
(a1),(a2),(b1),(b2),(c) and (d) by Remark \ref{rem:almost}.
Then it is easy to check that operators
\[
a_j':=\left(\frac{1}{1+\e}\right)^{1/2}a_j
\quad\text{for}\ j=1,\dots,p.
\]
satisfy the above properties for sufficiently small $\e>0$.
\erem


\blem\label{lem:almost}
Let $\xi$, $\eta$ be two almost $\delta$-related unit vectors 
in a Hilbert $M$-bimodule.
Then there exist $b_1,\dots, b_p\in M$ such that
\begin{enumerate}
\item[{\rm (1)}] 
${\dis |\ip{(1-\sum_{j=1}^pb_j^*b_j) \xi, \xi}|<4\delta^{1/2}}$;
\item[{\rm (2)}] 
${\dis |\ip{(1-\sum_{j=1}^pb_jb_j^*) \eta, \eta}|<4\delta^{1/2}}$;
\item[{\rm (3)}] 
${\dis \sum_{j=1}^p\|b_j\xi-\eta b_j\|^2<4\delta}$;
\item[{\rm (4)}] 
${\dis \sum_{j=1}^p\|b_j^*\eta-\xi b_j^*\|^2<4\delta}$.
\item[{\rm (5)}]
${\dis \sum_{j=1}^pb_j^*b_j\leq 1}$;
\item[{\rm (6)}]
${\dis \sum_{j=1}^pb_jb_j^*\leq 1}$.
\end{enumerate}
\elem

\bpf
Choose $0<\e<1$ such that
${\dis
5\e\leq 2\delta^{1/2}
}$
and 
${\dis
4\e^{1/2}\leq \delta^{1/2}.
}$
By Remark \ref{rem:almost2}, 
there exist $a_1,\dots,a_p$ of $M$ 
satisfying 
\begin{enumerate}
\item[(a'1)] ${\dis \|(1-\sum_{j=1}^pa_j^*a_j)\xi\|<\e^2,
\|(1-\sum_{j=1}^pa_j^*a_j)\eta\|<\e^2}$;
\item[(a'2)] ${\dis \|\xi(1-\sum_{j=1}^pa_j^*a_j)\|<\e^2,
\|\eta(1-\sum_{j=1}^pa_j^*a_j)\|<\e^2}$;
\item[(b'1)] ${\dis \|(1-\sum_{j=1}^pa_ja_j^*)\xi\|<\e^2,
\|(1-\sum_{j=1}^pa_ja_j^*)\eta\|<\e^2}$;
\item[(b'2)] ${\dis \|\xi(1-\sum_{j=1}^pa_ja_j^*)\|<\e^2,
\|\eta(1-\sum_{j=1}^pa_ja_j^*)\|<\e^2}$;
\item[(c')] ${\dis \sum_{j=1}^p\|a_j\xi-\eta a_j\|^2<\delta}$;
\item[(d')] ${\dis \sum_{j=1}^p\|a_j^*\eta-\xi a_j^*\|^2<2\delta}$;
\item[(e')] ${\dis \|(\sum_{j=1}^pa_j^*a_j)\xi\|\leq 1,
\|(\sum_{j=1}^pa_j^*a_j)\eta\|\leq 1}$;
\item[(f')] ${\dis \|(\sum_{j=1}^pa_ja_j^*)\xi\|\leq 1,
\|(\sum_{j=1}^pa_ja_j^*)\eta\|\leq 1}$.
\end{enumerate}
Then we define cp maps $S$ and $T$ on $M$ by
\[
S(x)\coloneqq\sum_{j=1}^pa_j^*xa_j,
\quad
T(x)\coloneqq\sum_{j=1}^pa_jxa_j^*
\quad
\text{for}\ x\in M.
\]
We define $e\coloneqq1_{[1-\e, 1+\e]}(S(1))$
and $f\coloneqq1_{[1-\e, 1+\e]}(T(1))$.
Since $(S(1)-1)^2\geq\e^2(1-e)$,
we have
\[
\e^2\|(1-e)\xi\|^2
=\e^2\ip{(1-e)\xi, \xi} 
\leq\|(S(1)-1)\xi\|^2<\e^4.
\]
Hence
\[
\|(1-e)\xi\|\leq\e.
\]
Similarly we have
\[
\|(1-e)\eta\|\leq\e,\quad
\|\xi(1-e)\|\leq\e,\quad
\|\eta(1-e)\|\leq\e.
\]
We also obtain 
\[
\|(1-f)\xi\|\leq\e,\quad
\|(1-f)\eta\|\leq\e,\quad
\|\xi(1-f)\|\leq\e,\quad
\|\eta(1-f)\|\leq\e.
\]
Next we define cp maps $S'$ and $T'$ on $M$ by
\[
S'(x)\coloneqq\frac{1}{1+\e}eS(x)e,
\quad
T'(x)\coloneqq\frac{1}{1+\e}fT(x)f,
\quad
\text{for}\ x\in M.
\]
Then $S'(1)\leq 1$ and $T'(1)\leq 1$.
In particular, $S'$ and $T'$ are contractive.

Now we define 
\[
b_j\coloneqq\frac{1}{\sqrt{1+\e}}fa_je,
\quad 
\text{for}\
1\leq j\leq p.
\]
Then
\[
\sum_{j=1}^pb_j^*b_j
=
\frac{1}{1+\e}\sum_{j=1}^pea_j^*fa_je
=
S'(f)
\leq
1.
\]
Similarly we have ${\dis \sum_{j=1}^pb_jb_j^*\leq 1}$.
Thus we obtain (5) and (6).

We will check (1). 
Since $\xi$ and $\eta$ are unit vectors,
we have
\begin{align*}
&|\ip{(1-\sum_{j=1}^pb_j^*b_j)\xi, \xi}|
=
\left|\ip{\eta, \eta}-\frac{1}{1+\e}\sum_{j=1}^p\ip{ea_j^*fa_je\xi, \xi}\right|
\\
&\quad\leq\e+
\left|
\ip{\eta (1-\sum_{j=1}^pa_ja_j^*), \eta}\right|
+
\left|\ip{(1-f)\eta, \eta \sum_{j=1}^pa_ja_j^*}\right|
+
\left|\sum_{j=1}^p\ip{f\eta a_j, (\eta a_j-a_j\xi)}\right|
\\
&\quad+
\left|\sum_{j=1}^p\ip{f(\eta a_j-a_j\xi), a_j\xi}\right|
+
\left|\sum_{j=1}^p\ip{fa_j\xi, a_j(1-e)\xi}\right|
+
\left|\sum_{j=1}^p\ip{fa_j(1-e)\xi, a_je\xi}\right|
\\
&\quad\leq\e+
\|\eta (1-\sum_{j=1}^pa_ja_j^*)\|\|\eta\|
+
\|(1-f)\eta\|\|\eta\sum_{j=1}^pa_ja_j^*\|
\\
&\quad+
\left(\sum_{j=1}^p\|\eta a_j\|^2\right)^{1/2}
\left(\sum_{j=1}^p\|\eta a_j-a_j\xi\|^2\right)^{1/2}
+
\left(\sum_{j=1}^p\|\eta a_j-a_j\xi\|^2\right)^{1/2}
\left(\sum_{j=1}^p\|a_j\xi\|^2\right)^{1/2}
\\
&\quad+
\left(\sum_{j=1}^p\|a_j\xi\|^2\right)^{1/2}
\left(\sum_{j=1}^p\|a_j(1-e)\xi\|^2\right)^{1/2}
+
\left(\sum_{j=1}^p\|a_j(1-e)\xi\|^2\right)^{1/2}
\left(\sum_{j=1}^p\|a_je\xi\|^2\right)^{1/2}.
\end{align*}
By (e'), (f'), we have
$\sum_{j=1}^p\|\eta a_j\|^2\leq 1$,
and
$\sum_{j=1}^p\|a_j\xi\|^2\leq 1$.
Since the projection $e$ commutes with ${\dis \sum_{j=1}^pa_j^*a_j}$,
we have
\[
\sum_{j=1}^p\|a_j(1-e)\xi\|^2
=
\ip{\sum_{j=1}^pa_j^*a_j(1-e)\xi, (1-e)\xi}
=
\ip{\sum_{j=1}^pa_j^*a_j\xi, (1-e)\xi}
\leq
\|(1-e)\xi\|\leq\e,
\]
and
\[
\sum_{j=1}^p\|a_je\xi\|^2
=
\ip{\sum_{j=1}^pa_j^*a_je\xi, e\xi}
=
\ip{\sum_{j=1}^pa_j^*a_j\xi, e\xi}
\leq 1.
\]
Therefore
\[
|\ip{(1-\sum_{j=1}^pb_j^*b_j)\xi, \xi}|
\leq
5\e+2\delta^{1/2}\leq4\delta^{1/2}.
\]
Similarly we have (2).

Next we will check (3).
\begin{align*}
&(\sum_{j=1}^p\|b_j\xi-\eta b_j\|^2)^{1/2}
\leq
(\sum_{j=1}^p\|fa_je\xi-\eta fa_je\|^2)^{1/2}
\\
&\quad=
(\sum_{j=1}^p\|fa_j(e-1)\xi\|^2)^{1/2}
+
(\sum_{j=1}^p\|fa_j\xi(1-e)\|^2)^{1/2}
+
(\sum_{j=1}^p\|f(a_j\xi-\eta a_j)e)\|^2)^{1/2}
\\
&\quad\quad+
(\sum_{j=1}^p\|(f-1)\eta a_je)\|^2)^{1/2}
+
(\sum_{j=1}^p\|\eta(1-f) a_je)\|^2)^{1/2}
\\
&\quad\leq
(\ip{S(f)(e-1)\xi, (e-1)\xi})^{1/2}
+
(\ip{S(f)\xi(1-e), \xi(1-e)})^{1/2}
\\
&\quad\quad+
(\sum_{j=1}^p\|a_j\xi-\eta a_j\|^2)^{1/2}
\\
&\quad\quad+
(\ip{(f-1)\eta, (f-1)\eta T(e)})^{1/2}
+
(\ip{\eta(1-f), \eta (1-f)T(e)})^{1/2}
\\
&\quad<
4\e^{1/2}+\delta^{1/2}\leq2\delta^{1/2}.
\end{align*}
Similarly we can check (4).
The proof is complete.
\epf
 
 Then we can show Theorem \ref{thm:almost} 
 by using following lemma, which is 
 also proved by the similar arguments as in \cite[Lemma 2.5]{ha2}.

\blem[cf. {\cite[Lemma 2.5]{ha2}}]\label{lemma 2.5}
Assume that $\delta>0$ and $r\in\N$
satisfy
\[
\delta^{1/2}<\frac{1}{8r}.
\]
Let $\xi$, $\eta$ be two almost $\delta$-related
unit vectors in a Hilbert $M$-bimodule.
Then  
there exist $r$ operators $c_1,\dots,c_r\in M$
such that $\|c_j\|\leq 1$, $1\leq j\leq r$
and 
\begin{align*}
&\left\|\left(\sum_{j=1}^rc_j^*c_j-1\right)\xi\right\|^2<\frac{12}{r},
\quad
\left\|\left(\sum_{j=1}^rc_jc_j^*-1\right)\eta\right\|^2<\frac{12}{r}, \\
&\sum_{j=1}^p\|c_j\xi-\eta c_j\|^2<32\delta,
\quad
\sum_{j=1}^p\|c_j^*\eta-\xi c_j^*\|^2<32\delta.
\end{align*}
\elem

By similar arguments as in \cite{ha2}
with Lemma \ref{lemma 2.5},
we can prove almost unitary equivalence 
for two almost $\delta$-related $n$-tuples
in a Hilbert bimodule,
which is a generalization of \cite[Theorem 2.3]{ha2}.

\bthm[cf.\ {\cite[Theorem 2.3]{ha2}}]\label{thm:almost}
For every $n\in\N$ and every $\e>0$,
there exists a $\delta=\delta(n, \e)>0$
such that for all von Neumann algebras $M$
and all almost $\delta$-related $n$-tuples
$(\xi_1,\dots,\xi_n)$, $(\eta_1,\dots,\eta_n)$
of unit vectors in a Hilbert $M$-bimodule,
there exists a unitary $u\in M$
such that
\[
\|u\xi_k-\eta_ku\|<\e
\quad\text{for}\ 1\leq k\leq n.
\]
\ethm


\section{$\Gamma$-stable states}


\bdf[cf.\ {\cite[Definition 4.1]{ha4}}]
Let $\Gamma$ be a multiplicative subgroup of $\R^+$.
A fn state $\varphi$ on a von Neumann algebra $M$
is called {\em $\Gamma$-stable}
if for every $n\in\N$, $0<r\leq 1$ 
and $\gamma_1,\dots,\gamma_n\in\Gamma$
with $1=r\gamma_1+\cdots+r\gamma_n$,
there exist $n$ partial isometries $v_1,\dots,v_n\in M$
and a projection $e\in M$
such that 
\[
\sum_{j=1}^nv_jv_j^*=1,
\quad\varphi(e)=r,
\]
and
\[
e=v_j^*v_j,\quad
\varphi v_j=\gamma_jv_j\varphi
\quad\text{for}\ 1\leq j\leq n.
\]
\edf

\brem
If $\Gamma=\Q^+$, then 
$\Q$-stable states in \cite[Definition 4.1]{ha4}
is equivalent to our $\Q^+$-stable states.
Indeed, a fn state $\varphi$ is $\Q$-stable
in the sense of \cite[Definition 4.1]{ha4}
if and only if 
for $q_1,\dots,q_n\in\Q^+$ with $1=q_1+\cdots+q_n$,
there exist $n$ isometries $v_1,\dots,v_n\in M$
such that
\[
\sum_{j=1}^nv_jv_j^*=1,
\quad\text{and}\quad
\varphi v_j=q_jv_j\varphi,
\]
because of \cite[Lemma 4.6]{ha4}.
Therefore if $\varphi$ is $\Q^+$-stable,
then for $q_1,\dots,q_n\in\Q^+$ 
with $1=q_1+\cdots+q_n$,
there are partial isometries $v_1,\dots,v_n\in M$
and a projection $e\in M$
such that 
\[
\sum_{j=1}^nv_jv_j^*=1,
\quad\varphi(e)=1,
\]
and
\[
e=v_j^*v_j,\quad
\varphi v_j=q_jv_j\varphi
\quad\text{for}\ 1\leq j\leq n.
\]
Since $\varphi(e)=1$ and $\varphi$ is faithful,
we have $e=1$. Hence $v_1,\dots,v_n$ are isometries,
and thus $\varphi$ is $\Q$-stable.

Conversely, let $0<r\leq 1$ 
and $\gamma_1,\dots,\gamma_n\in\Q^+$
with $1=r\gamma_1+\cdots+r\gamma_n$.
Then $r\in\Q^+$ and put 
$q_j\coloneqq r\gamma_j\in\Q^+$
for $1\leq j\leq n$.
By using $\Q$-stability,
there are isometries $v_1,\dots,v_n\in M$
such that
\[
\sum_{j=1}^nv_jv_j^*=1,
\quad\text{and}\quad
\varphi v_j=q_jv_j\varphi.
\]
Moreover there is an isometry $w$ such that 
$\varphi w=rw\varphi$.
We define partial isometries
$w_j\coloneqq v_jw^*$.
Then
$w_j^*w_j=ww^*$,
and
\[
\sum_{j=1}^nw_jw_j^*=\sum_{j=1}^nv_jv_j^*=1.
\]
Moreover
\[
\varphi w_j
=\varphi v_jw^*
=q_jr^{-1}v_jw^*\varphi
=\gamma_jw_j\varphi,
\]
and
\[
\varphi(w_j^*w_j)=\gamma_j^{-1}\varphi(w_jw_j^*)
=\gamma_j^{-1}\varphi(v_jv_j^*)
=\gamma_j^{-1}q_j\varphi(v_j^*v_j)=r.
\]
\erem

\blem[cf. {\cite[Theorem 4.5]{ha4}}]\label{lem:matrix}
Let $\varphi$ be a $\Gamma$-stable fn state on a von Neumann algebra $M$,
and let $0<r\leq 1$, $\gamma_1,\dots,\gamma_n\in\Gamma$
with $1=r\gamma_1+\cdots+r\gamma_n$.
Then there exists a type I$_n$ subfactor $F$ of $M$
such that  $\sigma_t^\varphi(F)=F$ for $t\in\R$ and
$\varphi|_F=\Tr_n(h\ \cdot\ )$,
where 
\[
h=\begin{bmatrix}
r\gamma_1 & & \\
& \ddots & \\
& & r\gamma_n
\end{bmatrix}.
\]
\elem

\bpf
There exist $n$ partial isometries $v_1,\dots,v_n\in M$
and a projection $e\in M$
such that 
\[
\sum_{j=1}^nv_jv_j^*=1,
\quad\varphi(e)=r,
\]
and
\[
e=v_j^*v_j,\quad
\varphi v_j=\gamma_jv_j\varphi
\quad\text{for}\ 1\leq j\leq n.
\]
Then $e_{jk}\coloneqq v_jv_k^*$ for $1\leq j,k\leq n$
give a system of matrix units.
Moreover we have
\[
\varphi(e_{jk})=\varphi(v_jv_k^*)
=
\gamma_j\varphi(v_k^*v_j)
=
\delta_{jk}r\gamma_j.
\]
Since
\[
\sigma_t^\varphi(e_{jk})
=
\sigma_t^\varphi(v_jv_k^*)
=
\gamma_j^{it}\gamma_k^{-it}v_jv_k^*
=
\gamma_j^{it}\gamma_k^{-it}e_{jk},
\]
we have $\sigma_t^\varphi(F)=F$ for $t\in\R$.
\epf

\brem\label{rem:F^c}
If $\tau$ is a tracial fn state on a type II$_1$ factor, 
then it is easy to see that $\tau$ is $\{1\}$-stable.

If $\varphi$ is a fn state on a type III$_\lambda$ factor
($0<\lambda<1$) for which $\sigma_{t_0}^\varphi=\id$,
where $t_0=-2\pi/\log\lambda$, 
then $\varphi$ is $\{\lambda^m\}_{m\in\Z}$-stable.
Indeed, let $0<r\leq 1$ 
and $\gamma_1,\dots,\gamma_n\in\Gamma$
with $1=r\gamma_1+\cdots+r\gamma_n$.
Put $\lambda_j\coloneqq r\gamma_j$.
Since $M_\varphi$ is a type II$_1$ factor,
we can choose a projection $e$
and mutually orthogonal projections $e_1,\dots, e_n$ in $M_\varphi$
with sum $1$
such that
\[
\varphi(e)=r,\quad\text{and}\quad
\varphi(e_j)=\lambda_j
\quad\text{for}\ 1\leq j\leq n.
\]
By using \cite[Lemma 4.2]{ha2},
there exist partial isometries 
$v_1,\dots,v_n$ in $M$
such that 
\[
e=v_j^*v_j,
\quad
e_j=v_jv_j^*,
\quad\text{and}\quad
\varphi v_j=\gamma_jv_j\varphi
\quad
\text{for}\
1\leq j\leq n.
\]

In these cases, by Lemma \ref{lem:matrix},
we obtain a finite dimensional subfactor $F$ of $M$
such that 
\[
\sigma_t^\varphi(F)=F
\quad\text{for}\ t\in\R.
\]
Note that it is equivalent to 
\[
\varphi=\varphi|_F\otimes\varphi|_{F^c},
\]
where $F^c\coloneqq F'\cap M$.
We expect that $\varphi|_{F^c}$ is also $\Gamma$-stable.
If $\tau$ is tracial,
then $\tau|_{F^c}$ is also tracial,
and thus is $\{1\}$-stable.
If $M$ is a type III$_\lambda$ factor
with a fn state $\varphi$ for which $\sigma_{t_0}^\varphi=\id$,
then $F^c$ is also a type III$_\lambda$ factor 
and $\sigma_t^{\varphi|_{F^c}}$ is the restriction 
of $\sigma_t^\varphi$ to $F^c$.
Hence $\varphi|_{F^c}$ is also $\{\lambda^m\}_{m\in\Z}$-stable.
In the case of type III$_1$ factors,
it is proved in \cite[Theorem 4.5]{ha4}
that if $\varphi$ is $\Q^+$-stable,
then $\varphi|_{F^c}$ is also $\Q^+$-stable.
\erem

\blem[cf. {\cite[Lemma 4.6]{ha4}}]\label{lem:unit}
If $\varphi$ is $\Gamma$-stable, 
then for $\gamma\in\Gamma$,
there exist $m\in\N$ and 
partial isometries $w_1,\dots, w_m\in M$
such that 
\[
\sum_{j=1}^mw_j^*w_j=1,
\quad\text{and}\quad
\varphi w_j=\gamma w_j\varphi
\quad\text{for}\ 1\leq j\leq m.
\]
\elem

\bpf
If $\gamma=1\in\Gamma$, then $m=1$ and $w_1=1$.
If $0<\gamma<1$, then set $r\coloneqq \gamma$.
By using $\Gamma$-stability for $1=r\gamma^{-1}$,
there exists a partial isometry $v\in M$
such that $vv^*=1$,
$\varphi v=\gamma^{-1} v\varphi$,
and $\varphi(v^*v)=\gamma$.
Then $w_1\coloneqq v^*$ satisfies the desired properties.
If $1<\gamma$, then there is $m\in\N$ such that
\[
m\gamma^{-1}\geq 1.
\]
Then take $0<r\leq 1$ such that
$r(m\gamma^{-1})=1$.
By using $\Gamma$-stability,
there exist partial isometries $v_1,\dots,v_m\in M$
such that 
\[
\sum_{j=1}^mv_jv_j^*=1,
\quad
\varphi v_j=\gamma^{-1}v_j\varphi
\quad\text{for}\ 1\leq j\leq m.
\]
Then $w_1\coloneqq v_1^*,\dots,w_m\coloneqq v_m^*$ satisfy the desired properties.
\epf


\section{Injective factors and ITPFI factors}


Throughout this section,
we assume that $\Gamma$ is a multiplicative subgroup of $\R^+$,
which is $\{1\}$, $\{\lambda^m\}_{m\in\Z}$ with $0<\lambda<1$
or $\Q^+$.
We also assume that $M$ is an injective factor $M$
of non-type I
with separable predual
and 
$\varphi$ is a $\Gamma$-stable fn state on $M$
with $B(M, \varphi)=\C1$.
Then if $M$ is of type II$_1$,
then $\varphi$ is tracial with $\Gamma=\{1\}$,
and if $M$ is of type III$_\lambda$ ($0<\lambda<1$),
then $\sigma_{t_0}^\varphi=\id$ with 
$\Gamma=\{\lambda^m\}_{m\in\Z}$,
where $t_0=-2\pi/\log\lambda$.
If $M$ is of type III$_1$,
then we assume that $\varphi$ is $\Q^+$-stable.
We recall that every type III$_1$ factor with separable predual has a $\Q^+$-stable fn state by \cite[Theorem 4.2]{ha4}.
Moreover every injective type III$_1$ with separable predual factor has trivial bicentralizer by \cite[Theorem 2.3]{ha3}.

We prove the main theorem in this section.

\bthm\label{thm:itpfi}
Let $M$ be an injective factor $M$
with separable predual
and 
$\varphi$ be a $\Gamma$-stable fn state on $M$
with $B(M, \varphi)=\C1$.
Then $M$ is ITPFI.
\ethm

\blem[cf. {\cite[Lemma 5.4]{ha4}}]\label{lem:lem1}
Let $\varphi$ be a $\Gamma$-stable fn state 
with $B(M, \varphi)=\C1$
on an injective factor $M$.
Let $u_1,\dots,u_n\in\cU(M)$ and $\delta>0$.
Then there exist a finite dimensional $\sigma^\varphi$-invariant
subfactor $F$ of $M$ and $v_1,\dots,v_n\in\cU(F)$
satisfying the following:
for $\e>0$, there exist
operators $b_1,\dots,b_p$ in $M$ such that
for $1\leq k\leq n$,
\begin{enumerate}
\renewcommand{\labelenumi}{\rm(\alph{enumi})}\renewcommand{\itemsep}{0pt}
\item ${\dis \|(1-\sum_{j=1}^pb_j^*b_j)u_k\xi_\varphi\|<\e}$,
${\dis \|(1-\sum_{j=1}^pb_j^*b_j)v_k\xi_\varphi\|<\e}$,
\item ${\dis \|(1-\sum_{j=1}^p\e_{F,\varphi}(b_jb_j^*))u_k\xi_\varphi\|<\e}$,
${\dis \|(1-\sum_{j=1}^p\e_{F,\varphi}(b_jb_j^*))v_k\xi_\varphi\|<\e}$,
\item ${\dis \sum_{j=1}^p\|b_j\xi_\varphi-\xi_\varphi b_j\|^2<\delta}$,
\item ${\dis \sum_{j=1}^p\|b_ju_k-v_kb_j\|_\varphi^2<\delta}$.
\end{enumerate}
\elem

\bpf
By Theorem \ref{thm:semidiscrete}
and Remark \ref{rem:semidiscrete},
for $1>\delta>0$
there exists a ucp map $T\colon\bM_m\to M$
and $v_1,\dots,v_n\in\cU(\bM_m)$
such that a fn state 
$\psi:=\varphi\circ T=\tr_m(h_\psi\ \cdot\ )$ 
on $\bM_m$ 
satisfies

\begin{align*}
&\|\sigma_t^\varphi\circ T-T\circ\sigma_t^\psi\|
\leq\delta|t|
\quad\text{for}\ t\in\R,
\\
&\|T(v_k)-u_k\|_\varphi<\frac{\delta^{1/2}}{2}
\quad \text{for}\ 1\leq k\leq n,
\end{align*}
and
\[
\lambda_1/\lambda_2\in\Gamma
\quad\text{for}\
\lambda_1, \lambda_2\in\Sp(h_\psi).
\]
Since $\varphi$ is $\Gamma$-stable,
as in the proof of \cite[Lemma 5.4]{ha4}, we may assume that 
$F:=\bM_m\subset M$, and $T\colon F\to M$
satisfies
$\varphi\circ T=\varphi|_F$ and 
\[
\|\sigma_t^\varphi\circ T-T\circ\sigma_t^{\varphi|_F}\|
\leq\delta|t|,
\quad t\in\R.
\]
Set $\xi_k:=u_k\xi_\varphi$ and $\eta_k:=v_k\xi_\varphi$
for $0\leq k\leq n$, 
where $u_0=v_0=1$.

By \cite[Proposition 2.1]{ha1},
there exist $a_1,\dots,a_p\in M$
such that
\[
T(x)=\sum_{j=1}^pa_j^*xa_j
\quad\text{for}\
x\in F.
\]
Since $T$ is unital,
we have 
\[
\sum_{j=1}^pa_j^*a_j=1.
\]

If $M$ is of type II$_1$,
then operators $a_1,\dots, a_p$
satisfies the desired properties
by the proof of \cite[Proposition 5.2]{ha1}.

Next we assume that $M$ is properly infinite. 
By \cite[Proposition 2.1]{ha1},
we can take a single operator $a\in M$ such that
\[
T(x)=a^*xa
\quad\text{for}\
x\in F.
\]
Then we can find finitely many operators $a_j$ satisfy the following conditions.
\begin{enumerate}
\renewcommand{\labelenumi}{\rm(\alph{enumi}')}\renewcommand{\itemsep}{0pt}
\item ${\dis \|(1-\sum_{j=-p}^p a_j^*a_j)\xi_k\|<\e/2}$,
${\dis \|(1-\sum_{j=1}^p a_j^*a_j)\eta_k\|<\e/2}$.
\item ${\dis \|(1-\sum_{j=-p}^p \lambda^{-j}\e_{F,\varphi}(a_ja_j^*))\xi_k\|<\e/2}$,
${\dis \|(1-\sum_{j=-p}^p \lambda^{-j}\e_{F,\varphi}(a_ja_j^*))\eta_k\|<\e/2}$.
\item ${\dis \sum_{j=-p}^p\|a_j\xi_\varphi-\lambda^{-j/2}\xi_\varphi a_j\|^2<\delta}$.
\item ${\dis \sum_{j=-p}^p\|a_ju_k-v_ka_j\|^2_\varphi<\delta}$.
\end{enumerate}
If $M$ is of type III$_1$ and $\varphi$ is $\Q_+$-stable,
then by the argument of \cite[Lemma 5.4]{ha4},
we can choose $\lambda\in\Q_+$
and operators $a_j$ for $-p\leq j\leq p$ 
satisfy the conditions (a')-(d').

If $M$ is of type III$_\lambda$ ($0<\lambda<1$), 
then it suffices to set $a_j\coloneqq \e_j(a)$,
where $\e_j$ is the projection of norm one of  $M$
onto
\[
M_j=\{x\in M\mid \sigma_t^\varphi(x)=\lambda^{ijt}(x),
t\in\R\}
\]
given by
\[
\e_j(x)\coloneqq\frac{1}{t_0}\int_0^{t_0}
\sigma_t^\varphi(x)\lambda^{-ijt}\,dt.
\]
Indeed, it follows from similar arguments as in \cite[Lemma 5.4]{ha4}.
We give a sketch proof below. 
Note that every $x\in M$ has a formal expansion
\[
x\sim\sum_{j=-\infty}^\infty\e_j(x).
\]
For $\xi\in L^2(M,\varphi)$, we have
\[
\sum_{j=-\infty}^\infty\|\e_j(a)\xi\|^2
=
\frac{1}{t_0}\int_0^{t_0}\|
\sigma_t^\varphi(a)\xi\|^2\,dt
=\|\xi\|^2.
\]
Let $x\in F$.
Since
\begin{align*}
\varphi(\e_{F,\varphi}(\e_j(a)\e_j(a)^*)x)
&=
\varphi\circ\e_{F,\varphi}(\e_j(a)\e_j(a)^*x)
\\
&=
\varphi(\e_j(a)\e_j(a)^*x)
\\
&=
\varphi(\e_j(a)^*x\sigma_{-i}^\varphi(\e_j(a)))
\\
&=
\lambda^j\varphi(\e_j(a)^*x\e_j(a)),
\end{align*}
we obtain
\begin{align*}
\sum_{j=-\infty}^\infty\lambda^{-j}
\varphi(\e_{F,\varphi}(\e_j(a)\e_j(a)^*)x)
&=
\sum_{j=-\infty}^\infty
\varphi(\e_j(a)^*x\e_j(a))
\\
&=
\sum_{j=-\infty}^\infty
\ip{x\e_j(a)\xi_\varphi, \e_j(a)\xi_\varphi}
\\
&=
\frac{1}{t_0}\int_0^{t_0}
\ip{x\sigma_t^\varphi(a)\xi_\varphi, \sigma_t^\varphi(a)\xi_\varphi}\,dt
\\
&=
\frac{1}{t_0}\int_0^{t_0}
\varphi(\sigma_t^\varphi\circ T\circ\sigma_{-t}^\varphi(x))\,dt
\\
&=\varphi(x).
\end{align*}
Since
$\xi_\varphi\e_j(a)=\lambda^{j/2}\e_j(a)\xi_\varphi$,
we have
\[
\sum_{j=-\infty}^\infty
\|\e_j(a)\xi_\varphi-\lambda^{-j/2}\xi_\varphi\e_j(a)\|^2
=0.
\]
For $\xi, \eta\in L^2(M,\varphi)$, we obtain
\begin{align*}
\sum_{j=-\infty}^\infty
\ip{x\e_j(a)\xi, \e_j(a)\eta}
&=
\frac{1}{t_0}\int_0^{t_0}
\ip{x\sigma_t^\varphi(a)\xi, \sigma_t^\varphi(a)\eta}\,dt
\\
&=
\frac{1}{t_0}\int_0^{t_0}
\ip{\sigma_t^\varphi\circ T\circ\sigma_{-t}^\varphi(x)\xi, \eta}\,dt.
\end{align*}
Hence
\[
\|T(x)-\sum_{j=-\infty}^\infty
\e_j(a)^*x\e_j(a)\|\leq\delta\|x\|.
\]
Therefore it follows from the above arguments that, for sufficiently large integer $p>0$,
operators $a_j$ for $-p\leq j\leq p$ satisfy the conditions
(a')-(d').

By Remark \ref{rem:F^c},
$\varphi|_{F^c}$ is also $\Gamma$-stable.
By Lemma \ref{lem:unit}, 
for each $-p\leq j\leq p$,
there exist a finite set of operators
$c_{j,1},\dots,c_{j,p(j)}$ in $F^c$ such that 
\[
\varphi c_{j,l}=\lambda^{-j}c_{j,l}\varphi
\quad\text{and}\quad
\sum_{l=1}^{p(j)}c_{j,l}^*c_{j,l}=1.
\]
Then operators $b_{j,l}:=c_{j,l}a_j$ for $-p\leq j\leq p$ and $l=1,\dots,p(j)$ satisfy the desired properties
as in the proof of \cite[Lemma 5.4]{ha4}.
\epf

\blem[cf. {\cite[Lemma 5.5]{ha4}}]\label{lem:lem3}
Let $\delta>0$ and $u_1,\dots,u_n\in\cU(M)$,
Then there exist a finite dimensional 
$\sigma^\varphi$-invariant subfactor $F$ of $M$
and unitaries $v_1,\dots, v_n\in\cU(F)$
such that $n+1$-tuples of unit vectors
$(\xi_\varphi,u_1\xi_\varphi,\dots,u_n\xi_\varphi)$
and $(\xi_\varphi,v_1\xi_\varphi,\dots,v_n\xi_\varphi)$
are almost $\delta$-related.
\elem

\bpf
The proof is the same as \cite[Lemma 5.5]{ha4}.
So we only give a sketch of the proof.
Let $\delta>0$ and $u_1,\dots,u_n\in\cU(M)$.
Thanks to Lemma \ref{lem:lem1},
there exist a finite dimensional $\sigma^\varphi$-invariant
subfactor $F$ of $M$ and $v_1,\dots,v_n\in\cU(F)$
as in the statement of Lemma \ref{lem:lem1}.
For a given $\e>0$, we can choose
operators $b_1,\dots,b_p$ in $M$ 
satisfying the conditions (a)-(d) in Lemma \ref{lem:lem1}.

Let $\delta'>0$ be arbitrary.
Set $b:=\sum_{j=1}^pb_jb_j^*$.
Since $B(M,\varphi)=\C$,
thanks to \cite{ha3} and \cite[Proposition 2.6]{ha4}
we have
\[
\e_{F,\varphi}(b)\in\overline{\Conv}\{wbw^*
\mid w\in\cU(F^c), \|w\xi_\varphi-\xi_\varphi w\|<\delta'\},
\]
where the closure means the $\sigma$-strong operator topology.
Hence there exist 
$w_1,\dots,w_q\in\cU(F^c)$
and
$\lambda_1,\dots,\lambda_q\in\R^+$ 
with $\sum_{l=1}^q\lambda_q=1$
such that
\[
\|w_l\xi_\varphi-\xi_\varphi w_l\|<\delta'
\]
and
\[
\|(1-\sum_{l=1}^q\lambda_lw_lbw_l^*)u_k\xi_\varphi\|<\e,
\quad\text{and}\quad
\|(1-\sum_{l=1}^q\lambda_lw_lbw_l^*)v_k\xi_\varphi\|<\e.
\]
Then one can easily check that operators
\[
a_{j,l}:=\lambda_l^{1/2}w_lb_j
\]
for $j=1,\dots,p$ and $l=1,\dots,q$
have the desired properties for the almost $\delta$-relatedness.
\epf


Now we prove our main theorem in this section.

\bpf[Proof of Theorem \ref{thm:itpfi}]
Let $u_1,\dots,u_n\in\cU(M)$,
and $\e>0$.
Then we take $\delta=\delta(n,\e/4)>0$
with properties in Theorem \ref{thm:almost}.
By Lemma \ref{lem:lem3}, we choose a finite dimensional 
$\sigma^\varphi$-invariant subfactor $F$ of $M$
and $v_1,\dots, v_n\in\cU(F)$
such that $n+1$-tuples of unit vectors
$(\xi_\varphi,u_1\xi_\varphi,\dots,u_n\xi_\varphi)$
and $(\xi_\varphi,v_1\xi_\varphi,\dots,v_n\xi_\varphi)$
are almost $\delta$-related.
Therefore by Theorem \ref{thm:almost},
there exists $w\in\cU(M)$
such that
\[
\|w\xi_\varphi-\xi_\varphi w\|<\frac{\e}{4},
\]
and
\[
\|wu_k\xi_\varphi-v_k\xi_\varphi w\|<\frac{\e}{4}
\quad\text{for}\
1\leq k\leq n.
\]
Then
\begin{align*}
\|u_k-w^*v_kw\|_\varphi
&=
\|w^*(wu_k-v_kw)\xi_\varphi\|
\\
&\leq
\|wu_k\xi_\varphi-v_k\xi_\varphi w\|
+\|v_k(\xi_\varphi w-w\xi_\varphi)\|
\\
&<\frac{\e}{2}.
\end{align*}
Set $F_0\coloneqq w^*Fw$ and 
$w_k\coloneqq w^*v_kw$
for $1\leq k\leq n$.
Then
\[
\|u_k-w_k\|_\varphi<\frac{\e}{2}
\quad\text{for}\
1\leq k\leq n.
\]
Put $\varphi_0\coloneqq w^*\varphi w$.
Then since $F$ is $\sigma^\varphi$-invariant,
$F_0$ is also $\sigma^{\varphi_0}$-invariant,
i.e.,
\[
\varphi_0=\varphi_0|_{F_0}\otimes\varphi_0|_{F_0^c}
\]
Since $\xi_{\varphi_0}=w^*\xi_\varphi w$,
we have
\begin{align*}
\|\varphi-\varphi_0\|
&\leq
\|\xi_\varphi-w^*\xi_\varphi w\|
\|\xi_\varphi+w^*\xi_\varphi w\|
\\
&\leq 
2\|w\xi_\varphi-\xi_\varphi w\|
\\
&<\e.
\end{align*}
Therefore by \cite[Lemma 7.6]{cw},
$M$ is ITPFI.
\epf



\begin{thebibliography}{Was}


\bibitem[AH]{ah} H. Ando, U. Haagerup;
{\it Ultraproducts of von Neumann algebras.} 
J. Funct. Anal. {\bf 266} (2014), no. 12, 6842--6913. 

\bibitem[AW]{aw} H. Araki, E. J. Woods; 
{\it A classification of factors.} 
Publ. Res. Inst. Math. Sci. Ser. A. {\bf 4} (1968), 51--130.

\bibitem[Co1]{c3} A. Connes; 
{\it Classification of injective factors.} 
Ann. Math. {\bf 104} (1976), 73--115. 

\bibitem[Co2]{c4} A. Connes; 
{\it Factors of type III$_1$, property $L_\lambda'$, and closure of inner automorphisms.} 
J. Operator Theory {\bf 14} (1985), 189--211. 

\bibitem[CS]{cs} A. Connes, E. St\o rmer;
{\it Homogeneity of the state space of factors of type III$_1$.} 
J. Funct. Anal. {\bf 28} (1978), 187--196. 

\bibitem[CT]{ct} A. Connes, M. Takesaki; 
{\it The flow of weights on factors of type III.} 
Tohoku Math. J. {\bf 29} (1977), 473--575. 

\bibitem[CW]{cw} A. Connes, E. J. Woods; 
{\it Approximately transitive flows and ITPFI factors.} 
Ergod. Th. Dynam. Sys. {\bf 5} (1985), 203--236.

\bibitem[Ha1]{ha1} U. Haagerup; 
{\it A new proof of the equivalence of injectivity and hyperfiniteness for factors on a separable Hilbert space.}
J. Funct. Anal. {\bf 62} (1985), no. 2, 160--201.

\bibitem[Ha2]{ha2} U. Haagerup; 
{\it The injective factors of type III$_\lambda$, $0<\lambda<1$.}
Pacific. J. Math. {\bf 137} (1989), 265--310.

\bibitem [Ha3]{ha3} U. Haagerup;
{\it Connes' bicentralizer problem and uniqueness of the injective factor of type III$_1$.} 
Acta Math. {\bf 158} (1987), no. 1-2, 95--148. 

\bibitem[Ha4]{ha4} U. Haagerup; 
{\it On the uniqueness of the injective III$_1$ factor.} 
Doc. Math. {\bf 21} (2016), 1193--1226.

\bibitem[HS]{hs}  U. Haagerup, E. St\o rmer;
{\it Equivalence of normal states on von Neumann algebras and the flow of weights.} 
Adv. Math. {\bf 83} (1990), no. 2, 180--262.
 
\bibitem[HI]{hi} C. Houdayer, Y. Isono;
{\it Unique prime factorization and bicentralizer problem for a class of type III factors.} 
Adv. Math. {\bf 305} (2017), 402--455. 

\bibitem[Oc]{oc} A. Ocneanu;
{\it Actions of discrete amenable groups on von Neumann algebras.}
Lecture Notes in Mathematics {\bf 1138} 
Springer-Verlag, Berlin, 1985. iv+115 pp. 

\bibitem[Po]{po} S. Popa;
{\it A short proof of ``injectivity implies hyperfiniteness'' for finite von Neumann algebras.}
J. Operator Theory {\bf 16} (1986), no. 2, 261--272.

\bibitem[Tak]{tak} M. Takesaki;
{\it Theory of operator algebras. II.} 
Encyclopaedia of Mathematical Sciences, {\bf 125}. Operator Algebras and Non-commutative Geometry, {\bf 6}. Springer-Verlag, Berlin, 2003. xxii+518 pp.

\bibitem[Was]{was} S. Wassermann;
{\it Injective W$^*$-algebras.} 
Math. Proc. Cambridge Philos. Soc. {\bf 82} (1977), no. 1, 39--47.


\end{thebibliography}
\end{document}